\documentclass[12pt]{amsart}
\usepackage{amscd}
\usepackage{amssymb}
\usepackage{graphicx}
\usepackage{rotate,epic,eepic,epsf,epsfig}

\allowdisplaybreaks[1]

\newtheorem{thm}{Theorem}[section]
\newtheorem{rem}[thm]{Remark}

\newtheorem{defns}[thm]{Definitions}

\newtheorem{prop}[thm]{Proposition}

\newtheorem{defn}[thm]{Definition}
\theoremstyle{definition}
\newtheorem{exple}[thm]{Example}

\newtheorem{exples}[thm]{Examples}

\numberwithin{equation}{section}

\DeclareFixedFont{\petitefonte}{\encodingdefault}%
{\familydefault}{\seriesdefault}{\shapedefault}{5pt}

\newcommand{\toutpetit}{\tiny \smaller\smaller\smaller\smaller\smaller\smaller\smaller\smaller\smaller\smaller
 \smaller\smaller\smaller\smaller\smaller\smaller}

\newcommand{\C}{\mathbb C}

\newcommand{\n}{\noindent}
\renewcommand{\a}{\alpha}

\newcommand{\resumename}{R\'esum\'e}



\thanks{This work was supported by the CMCU contract 06/S 1502 (Partenariat Hubert Curien Utique),
O. Khlifi thanks the Universit\'e de Bourgogne and D. Arnal thanks the Facult\'e des Sciences de Sfax for
their kind hospitalities they received during their visits.}
\begin{document}
\title[Diamond representations for rank two semisimple Lie algebras]{Diamond representations of rank two semisimple Lie algebras}

\author[B. Agrebaoui, D. Arnal and O. Khlifi]{Boujemaa Agrebaoui, Didier Arnal and Olfa Khlifi}

\address{D\'epartement de Mathematiques, Facult\'e des Sciences de Sfax, Route de Soukra, km 3,5,
B.P. 1171, 3000 Sfax, Tunisie.}
\email{bagrebaoui@yahoo.fr}

\address{Institut de Math\'ematique de Bourgogne, UMR Cnrs 5584,
Universit\'e de Bourgogne, U. F. R. Sciences et Techniques\\
B. P. 47870, F- 21078 Dijon Cedex, France.}
\email{Didier.Arnal@u-bourgogne.fr}

\address{D\'epartement de Mathematiques, Facult\'e des Sciences de Sfax, Route de Soukra, km 3,5,
B.P. 1171, 3000 Sfax, Tunisie.}
\email{khlifi\_olfa@yahoo.fr}

\begin{abstract}

The present work is a part of a larger program to construct explicit combinatorial models for the
(indecomposable) regular representation of the nilpotent factor $N$ in the Iwasawa decomposition of a
semi-simple Lie algebra $\mathfrak g$, using the restrictions to $N$ of the simple finite dimensional
modules of $\mathfrak g$.\\

Such a description is given in \cite{[ABW]}, for the cas $\mathfrak g=\mathfrak{sl}(n)$. Here, we give
the analog for the rank 2 semi simple Lie algebras (of type $A_1\times A_1$, $A_2$, $C_2$ and $G_2$).
The algebra $\mathbb C[N]$ of polynomial functions on $N$ is a quotient, called reduced shape algebra
of the shape algebra for $\mathfrak g$. Basis for the shape algebra are known, for instance the so called
semi standard Young tableaux (see \cite{[ADLMPPrW]}). We select among the semi standard tableaux, the so
called quasi standard ones which define a kind basis for the reduced shape algebra.\\
\end{abstract}

\maketitle

\vspace{.20cm}

\section{Introduction}

\

We will study the diamond cone of representations for the nilpotent factor $N^+$ of any rank 2 semi simple
Lie algebra $\mathfrak g$. This is the indecomposable regular representation onto $\mathbb C[N^-]$, described
from explicit realizations of the restrictions to $N^-$ of the simple $\mathfrak g$-modules $V^\lambda$.\\

In \cite{[ABW]}, this description is explicitely given in the case $\mathfrak g=\mathfrak{sl}(n)$, using the
notion of quasi standard Young tableaux. Roughly speaking, a quasi standard Young tableau is an usual semi
standard Young tableau such that, it is impossible to extract the top of the first column, either because
this top of column is not `trivial', {\sl i.e.} it does not consist of numbers $1,2,\dots,k$, or because,
when we push out this top by pushing the $k$ fist rows of the tableau, we do not get a semi standard tableau.\\

Let us come back for the case of rank 2 Lie algebra $\mathfrak g$. The modules $V^\lambda$ have well known
explicit realizations (see for instance \cite{[FH]}). They are characterized by their highest weight
$\lambda=a\omega_1+b\omega_2$. In \cite{[ADLMPPrW]}, there is a construction for a basis for each
$V^\lambda$, as the collection of all semi standard tableaux with shape $(a,b)$. The definition
and construction of semi standard tableaux for $\mathfrak g$ uses the notion of grid poset and their ideals.
It is possible to perform compositions of grid posets, the ideals of these compositions (of $a$ grid posets
associated to $V^{\omega_1}$ and $b$ grid posets associated to $V^{\omega_2}$) give a basis for $V^\lambda$
if $\lambda=a\omega_1+b\omega_2$.\\

Thus, we realize the Lie algebra $\mathfrak g$ as a subalgebra of $\mathfrak{sl}(n)$ (with $n=4,3,4,7$),
and we recall the notion of shape algebra for $\mathfrak g$, it is the direct sum of all the $V^\lambda$,
but we can see it as the algebra $\mathbb C[G]^{N^+}$ of all the polynomial functions on the group $G$
corresponding to $\mathfrak g$, which are invariant under right action by $N^+$. This gives a very concrete
interpretaion of the semi standard tableaux for $\mathfrak g$.\\

The algebra $\mathbb C[N^-]$ is the restriction to $N^-$ of the functions in $\mathbb C[G]$. But it is
also a quotient of the shape algebra by the ideal generated by $\renewcommand{\arraystretch}{0.7}\begin{array}{l}\framebox{$1$}\\
\framebox{2}\end{array}-1$, $\framebox{1}-1$. We call this quotient the reduced
shape algebra for $\mathfrak g$. To give a basis for this quotient, we define, case by case, the quasi
standard tableaux for $\mathfrak g$. They are semi standard Young tableaux, with an extra condition,
which is very similar to the condition given in the $\mathfrak{sl}(n)$ case. We prove the quasi standard
Young tableaux give a kind basis for the reduced shape algebra.\\

\maketitle

\section{Semi standard Young tableaux - quasi standard Young tableaux of $SL(n)$}

\vskip0.7cm

\subsection{Semi standard Young tableaux}

\n

Recall that the Lie algebra $\mathfrak{sl}(n) = \mathfrak{sl}(n,
\mathbb C)$ is the set of $n\times n$ traceless matrices, it is the
Lie algebra of the Lie group $SL(n)$ of $n\times n$ matrices, with
determinant 1.

\

Denote $N^+$ the subgroup of all the matrices $n^+ = \left(\begin{matrix}1&&&*\\
&& \ddots &&\\
0&&&1
\end{matrix}\right).$ Let us consider the algebra ${\mathbb
C}[SL(n)]^{N^+}$ of polynomial functions on the group $SL(n)$, which
are invariant under the right multiplication by the subgroup
$N^+.$
\begin{exple}

\n

Let $k<n$ and $1\leq i_1<i_2<\dots<i_k\leq n$. We define :

$$
\begin{array}{cccc}
\begin{array}{|c|}
\hline
i_1 \\
\hline
i_2 \\
\hline
\vdots \\
\hline
i_k \\
\hline
\end{array}:& SL(n)  &  \longrightarrow &
\mathbb{C}  \\
 & g   & \longmapsto   &
det(submatrix(g,(i_1...i_k,1...k))) \\
\end{array}$$

\n
{\sl i.e.} for an element $g \in SL(n),$ we associate the polynomial
function which is the determinant of the submatrix of $g$ obtained by considering
the $k$ first columns of $g$ and the rows $i_1,\dots,i_k$.

\end{exple}

If $k$ is fixed, $SL(n)$ acts on the vector space of all columns  $\begin{array}{|c|} \hline
i_1 \\
\hline
i_2 \\
\hline
\vdots \\
\hline
i_k \\
\hline
\end{array}$ as on $\wedge^{k} \mathbb{C}^{n}$.

Thus we look for $Sym^{\bullet}(\bigwedge\C^n)=Sym^{\bullet}
(\mathbb{C}^{n} \oplus \wedge^{2} \mathbb{C}^{n} \oplus \dots \oplus \wedge^{n-1} \mathbb{C}^{n}
)$. A basis for this algebra is given by the Young tableaux
$$
\begin{tabular}{|c|cccc}\hline
                $ i_{1}^{1}$ &   \multicolumn{1}{|c|}{$i_{1}^{2}$} & & $\cdots$  &
                   \multicolumn{1}{|c|}{$i_{1}^{r}$} \\  \hline
                   $\vdots$ &     \multicolumn{1}{|c|}{$\vdots$} &    &  \multicolumn{1}{c|}{$\vdots$}  &             \\     \cline{1-4}
    &   \multicolumn{1}{|c|}{$ i_{k_2}^{2}$} &        &           &\\     \cline{1-2}
 $ i_{k_1}^{1}$   &      &     &    &     \\   \cline{1-1}
\end{tabular}
$$
such that
~~ $k_1\geq k_2\geq ...\geq k_r$ and if $k_j= k_{j+1}$ then $\left(
\begin{array}{c}
 i _{1}^{j} \\
 \vdots  \\
 i_{k_j}^{j} \\
 \end{array}
\right) \leq  \left(\begin{array}{c}
 i _{1}^{j+1} \\
 \vdots  \\
 i_{k_j}^{j+1} \\
 \end{array}
\right)$ for the lexicographic ordering.

Recall now that the fundamental representations of $\mathfrak{sl}(n)$ are
the natural ones on ${\mathbb C}^n,~\wedge^2{\mathbb
C}^n,\dots,~\wedge^{n-1}{\mathbb C}^n$ and they have highest weights
$\omega_1,\dots,~\omega_{n-1}$.

From a Borel-Weyl theorem, we obtain that  each simple $\mathfrak{sl}(n)$-module has a highest
weight $\lambda$, there are non negative integral numbers $a_1,\dots,a_{n-1}$ such that
$$
\lambda=a_1\omega_1+\dots+a_{n-1}\omega_{n-1}
$$
and this highest weight characterizes the module.
Note ${\mathbb S}^\lambda$ (or $\Gamma_{a_1,\dots,a_{n-1}}$) this module,
it is a submodule of the tensor product
$$
Sym^{a_1}({\mathbb C}^n)\otimes Sym^{a_2}(\wedge^2{\mathbb C}^n)\otimes
\dots\otimes Sym^{a_{n-1}}(\wedge^{n-1}{\mathbb C}^n).
$$

The direct sum ${\mathbb S}^{\bullet}$ of all the simple modules $\mathbb S^\lambda$ is the shape algebra of
$SL(n)$. As an algebra, it is isomorphic to ${\mathbb C}[SL(n)]^{N^+}$ ( see \cite{[FH]}).\\

Now, we have a natural mapping from  $Sym^{\bullet}
(\mathbb{C}^{n} \oplus \wedge^{2} \mathbb{C}^{n} \oplus \dots \oplus \wedge^{n-1} \mathbb{C}^{n}
)$ onto  ${\mathbb C}[SL(n)]^{N^+}$ which is just the evaluation map:
$$
\begin{tabular}{|c|cccc}\hline
                $ i_{1}^{1}$ &   \multicolumn{1}{|c|}{$i_{1}^{2}$} & & $\cdots$  &
                   \multicolumn{1}{|c|}{$i_{1}^{r}$} \\  \hline
                   $\vdots$ &     \multicolumn{1}{|c|}{$\vdots$} &    &  \multicolumn{1}{c|}{$\vdots$}  &             \\     \cline{1-4}
    &   \multicolumn{1}{|c|}{$ i_{k_2}^{2}$} &        &           &\\     \cline{1-2}
 $ i_{k_1}^{1}$   &      &     &    &     \\   \cline{1-1}
\end{tabular} \longmapsto \left( g \longmapsto \begin{array}{|c|} \hline
i_1^1 \\
\hline
\vdots \\
\hline
 \\
\hline
i_{k_1}^1 \\
\hline
\end{array} ~(g) \cdot \begin{array}{|c|} \hline
i_1^2 \\
\hline
\vdots \\
\hline
i_{k_2}^2 \\
\hline
\end{array}~(g) \ldots \begin{array}{|c|} \hline
i_1^r \\
\hline \end{array}~(g) \right).
$$

Thanks to the Gauss method, all the $N^+$ right
invariant polynomial functions on $SL(n)$ are polynomial functions
in the polynomials $\begin{array}{|c|} \hline
i_1 \\
\hline
i_2 \\
\hline
\vdots \\
\hline
i_k \\
\hline
\end{array}$, thus :\\

\begin{prop}

 \n

The map from $Sym^{\bullet}(\bigwedge\mathbb{C}^n)=Sym^{\bullet}(\mathbb{C}^n\oplus\wedge^2
\mathbb{C}^n\oplus \dots \oplus \wedge^{n-1} \mathbb{C}^n)$ onto
$\mathbb S^\bullet={\mathbb C}[SL(n)]^{N^+}$ is a surjective mapping.
\end{prop}
\begin{defn}

\n

A Young tableaux of shape $\lambda$ is {\bf{semi standard}} if its entries are
increasing along each row (and strictly increasing along each
column).
\end{defn}
\begin{defn}

\n

Let $T$ be a tableau. If $T$ contains $a_i$ columns with height $i$
$(i=1,...,n-1),$ we call shape of $T$ the $(n-1)$-uplet
$\lambda(T)=(a_1,...,a_{n-1})$. We consider the patial ordering on the family of shapes $\mu=(b_1,\dots,b_{n-1})\leq\lambda=(a_1,\dots,a_{n-1})$ if and only if 
$$
b_1\leq a_1,~~\dots,~~b_{n-1}\leq a_{n-1}.
$$
\end{defn}

\

\begin{thm}

\n

\begin{itemize}
\item[1)] The algebra ${\mathbb S}^\bullet=
\bigoplus_{\lambda}~~{\mathbb S}^\lambda,$ is isomorphic to the
quotient of $Sym^{\bullet}(\bigwedge\mathbb{C}^n)$ by the kernel
$\mathcal{PL}$ of the evaluation mapping. This ideal is
generated by the Pl\"ucker relations.\\
\item[2)] If $\lambda=a_1\omega_1+\dots+a_{n-1}\omega_{n-1}$, a basis for ${\mathbb S}^\lambda$ is
given by the set of semi standard Young tableaux $T$ of shape
$\lambda$.
\end{itemize}
\end{thm}
\begin{exple} The $\mathfrak{sl}(3)$ case

\n

We have one and only one Pl\"ucker relation: $$
\renewcommand{\arraystretch}{0.7}{\begin{array}{l}
\framebox{1}\framebox{3}\\
\framebox{2} \\
\end{array}}- \renewcommand{\arraystretch}{0.7}{\begin{array}{l}
\framebox{2}\framebox{1}\\
\framebox{3} \\
\end{array}} + \renewcommand{\arraystretch}{0.7}{\begin{array}{l}
\framebox{1}\framebox{2}\\
\framebox{3} \\
\end{array}}=0.$$

\n Then to obtain a basis for the algebra ${\mathbb S}^\bullet$, we
reject exactly the non semi standard Young tableaux : the tableaux which contain
$\renewcommand{\arraystretch}{0.7}{\begin{array}{l}
\framebox{2}\framebox{1}\\
\framebox{3} \\
\end{array}}$ as a subtableau.
\end{exple}

The Cartan subalgebra $\mathfrak h$ of $ \mathfrak{sl}(n)$ is the $(n-1)$ dimensional vector space consisting of diagonal,
traceless matrices $H=(h_{ij}).$ The usual basis $(\alpha_1,\dots,\alpha_{n-1})$ of $\mathfrak h^\ast$
is given by simples roots
$\alpha_i= \lambda_i- \lambda_{i+1} ~~\hbox{where} ~~\lambda_i(H)=h_{ii}.$\\
$\mathfrak h^\ast$ is an Euclidean vector space with a scalar product given by the Killing form.
We shall draw pictures in the real vector space $\mathfrak h^\ast_{\mathbb{R}}$ generated by $\alpha_i$.

For $\mathfrak{sl}(3),$ we note $\alpha=\alpha_1$ and $\beta=\alpha_2$.\\

Following \cite{[ABW]}, we look at the action of the nilpotent group $^t N^+$ onto the highest weight vector
$v_\lambda$ in $\mathbb S^\lambda$. This action generates the representation space $\mathbb{S} ^{\lambda}$.
Thus, as basis for the dual of the Cartan
subalgebra, we choose the simple negative roots. The action of $X_{-\alpha}$ on a weight vector is pictured
by an arrow
\begin{picture}(50,20)(150,120)\toutpetit

\path(150,120)(190,120)

\put(170,120){\vector(1,0){1}}

\put(167,114){$\alpha$}
\put(191,120){.}
\end{picture}
\\
\begin{exple}

\n

With the convention above, we get the following weight diagrams of $\Gamma_{a,b}$
 for $\mathfrak{sl}(3)$, for $a+b \leq 2$ :
\vfill\eject

\n

\begin{center} \begin{tabular}{|c|c|c|}\hline

\begin{picture}(130,-20)(130,20)

 \put(190,13){$\Gamma_{0,0}:$}
\end{picture}            & \hskip2.5cm
\begin{picture}(50,-50)(200,60)\toutpetit
\put(190,58){\circle{3}} \put(188,50){$0$}   \end{picture}
\hskip2.5cm \\
&\\ \hline

\begin{picture}(130,50)(130,120)

 \put(190,140){$\Gamma_{1,0}:$}

\end{picture}
& \hskip2.5cm \begin{picture}(180,40)(160,117) \toutpetit

\path(195,120)(235,120)

\path(195,120)(215,150)

\put(216,120){\vector(1,0){1}}

\put(212,115){$\alpha$}

\put(235,120){\circle{3}}

\put(195,120){\circle{3}}

\put(181,118) {\framebox{1}}

\put(238,118){\framebox{2}}

\path(235,120)(215,150)

\put(224,136){\vector(-1,1){1}}

\put(229,135){$\beta$}

\put(215,150){\circle{3}}

\put(210,157){\framebox{3}}
\end{picture} \hskip2.5cm \\ \hline
\begin{picture}(130,50)(130,140)

 \put(190,170){$\Gamma_{0,1}:$}

\end{picture} & \hskip2.5cm

\begin{picture}(200,70)(120,65) \toutpetit

\path(170,120)(205,120)

\put(190,120){\vector(1,0){1}}

\put(185,123){$\alpha$}

\put(205,120){\circle{3}}

\put(170,120){\circle{3}}

\put(151,118) {$\begin{array}{l}
\framebox{1}\\
\framebox{3} \\
\end{array}$}

\put(203,118){$\begin{array}{l}
\framebox{2}\\
\framebox{3} \\
\end{array}$}

\path(170,120)(189,88)

\path(205,120)(189,88)

\put(189,89){\circle{3}}

\put(181,102){\vector(1,-1){1}}

\put(175,97){$\beta$}

\put(179,75){$\begin{array}{l}
\framebox{1}\\
\framebox{2} \\
\end{array}$}
\end{picture} \hskip2.5cm
\\   \hline
\begin{picture}(130,50)(130,120)

\put(190,155){$\Gamma_{2,0}:$}

\end{picture}& \hskip2.5cm

\begin{picture}(200,100)(130,72) \toutpetit

\path(170,120)(205,120)

\put(190,120){\vector(1,0){1}}

\put(185,124){$\alpha$}

\put(205,120){\circle{3}}

\put(170,120){\circle{3}}

\put(141,121) {$
\begin{array}{l}
\framebox{1}\framebox{3}
\end{array}$}

\put(204,121){$\begin{array}{l}
\framebox{2}\framebox{3}\\
\end{array}$}

\path(170,120)(189,88)

\path(205,120)(189,88)

\put(189,88){\circle{3}}

\put(180,102){\vector(-1,1){1}}

\put(175,97){$\beta$}

\put(176,78){$\begin{array}{l}

\framebox{1}\framebox{2}
\end{array}$}

\put(171,160){$\begin{array}{l}

\framebox{3}\framebox{3}
\end{array}$}

\put(222,88){$\begin{array}{l}

\framebox{2}\framebox{2}
\end{array}$}

\put(126,88){$\begin{array}{l}

\framebox{1}\framebox{1}
\end{array}$}
\path(154,88)(224,88)

\path(154,88)(186,152)

\path(186,152)(224,88)

\put(154,88){\circle{3}}

\put(224,88){\circle{3}}

\put(186,152){\circle{3}}

\put(216,88){\vector(1,0){1}}

\put(211,83){$\alpha$}

\put(170,88){\vector(1,0){1}}

\put(167,83){$\alpha$}

\put(215,102){\vector(-1,1){1}}

\put(220,101){$\beta$}

\put(195,137){\vector(-1,1){1}}

\put(199,137){$\beta$}
\end{picture} \hskip2.5cm \\  \hline
&  \\
&  \\
 \begin{picture}(130,50)(130,120)

\put(190,120){$\Gamma_{0,2}:$}

\end{picture} & \hskip2.5cm

\begin{picture}(180,50)(160,102) \toutpetit

\path(195,120)(235,120)

\path(195,120)(215,150)

\put(216,120){\vector(1,0){1}}

\put(212,115){$\alpha$}

\put(235,120){\circle{3}}

\put(195,120){\circle{3}}

\put(165,113) {$\begin{array}{l}
\framebox{1}\framebox{1}\\
\framebox{2}\framebox{3}\\
\end{array}$}

\put(235,116){$\begin{array}{l}
\framebox{1}\framebox{2}\\
\framebox{2}\framebox{3}\\
\end{array}$}
\path(235,120)(215,150)

\path(255,150)(215,150)

\path(175,150)(215,150)

\path(215,90)(255,150)

\path(175,150)(215,90)

\put(175,150){\circle{3}}

\put(255,150){\circle{3}}

\put(215,90){\circle{3}}

\put(224,136){\vector(-1,1){1}}

\put(200,100){$\beta$}

\put(205,105){\vector(-1,1){1}}

\put(185,135){\vector(-1,1){1}}

\put(177,133){$\beta$}

\put(195,150){\vector(1,0){1}}

\put(194,153){$\alpha$}

\put(235,150){\vector(1,0){1}}

\put(234,153){$\alpha$}

\put(230,135){$\beta$}

\put(215,150){\circle{3}}

\put(200,162){$\begin{array}{l}
\framebox{1}\framebox{2}\\
\framebox{3}\framebox{3}\\
\end{array}$}

\put(146,152){$\begin{array}{l}
\framebox{1}\framebox{1}\\
\framebox{3}\framebox{3}\\
\end{array}$}

\put(253,152){$\begin{array}{l}
\framebox{2}\framebox{2}\\
\framebox{3}\framebox{3}\\
\end{array}$}

\put(200,76){$\begin{array}{l}
\framebox{1}\framebox{1}\\
\framebox{2}\framebox{2}\\
\end{array}$}
\end{picture} \hskip2.5cm\\
&  \\
&  \\
&  \\  \hline
$ \begin{picture}(130,50)(130,170)

\put(190,150){$\Gamma_{1,1}:$}

\end{picture}  $  & \hskip2.5cm
\begin{picture}(180,50)(230,130) \toutpetit

\path(255,150)(335,150) \path(335,150)(375,90) \path(375,90)(335,30)
\path(255,30)(335,30)   \path(255,30)(215,90) \path(255,150)(215,90)
\path(295,90)(335,30) \path(295,90)(255,150)   \path(215,90)(375,90)
\path(255,30)(335,150)

\put(335,150){\circle{3}} \put(255,30){\circle{3}}
\put(215,90){\circle{3}} \put(335,30){\circle{3}}
\put(255,150){\circle{3}} \put(375,90){\circle{3}}
\put(295,90){\circle{3}}

\put(300,30){\vector(1,0){1}} \put(298,26){$\alpha$}

\put(231,65){\vector(-1,1){1}} \put(228,59){$\beta$}

\put(323,48){\vector(-1,1){1}} \put(325,49){$\beta$}

\put(268,130){\vector(-1,1){1}} \put(271,132){$\beta$}

\put(260,90){\vector(1,0){1}} \put(258,85){$\alpha$}

\put(330,90){\vector(1,0){1}} \put(328,85){$\alpha$}

\put(300,150){\vector(1,0){1}} \put(300,153){$\alpha$}

\put(351,125){\vector(-1,1){1}} \put(354,127){$\beta$}

\put(186,90) {$\begin{array}{l}
\framebox{1} \framebox{1}\\
\framebox{3}\\
\end{array}$}

\put(270,76){$\begin{array}{l}
\framebox{1}\framebox{3}\\
\framebox{2}\\
\end{array}$}

\put(234,155){$\begin{array}{l}
\framebox{1}\framebox{3}\\
\framebox{3}\\
\end{array}$}

\put(336,151){$\begin{array}{l}
\framebox{2}\framebox{3}\\
\framebox{3}\\
\end{array}$}

\put(334,20){$\begin{array}{l}
\framebox{1}\framebox{2}\\
\framebox{2}\\
\end{array}$}

\put(290,101){$\begin{array}{l}
\framebox{1}\framebox{2}\\
\framebox{3}\\
\end{array}$}

\put(373,87){$\begin{array}{l}
\framebox{2}\framebox{2}\\
\framebox{3}\\
\end{array}$}

\put(237,15){$\begin{array}{l}
\framebox{1}\framebox{1}\\
\framebox{2}\\
\end{array}$}
\hskip1.5cm
\end{picture} \hskip2.5cm \\
&  \\
&  \\
&  \\
&  \\
&  \\
&  \\
&  \\
&  \\
&  \\
&  \\  \hline

\end{tabular}
\end{center}
\end{exple}
\vfill\eject
\subsection{Quasi standard Young tableaux for $\mathfrak{sl}(n)$}

\n

We now are interested by the restriction of polynomial functions on $SL(n)$
to the subgroup $N^-=~^tN^+$. This restriction leads to an exact sequence: ( see \cite{[ABW]} )
$$
0 \longrightarrow \big< ~~\begin{array}{|c|} \hline
1 \\
\hline
2 \\
\hline
\vdots \\
\hline
k \\
\hline
\end{array}-1,~k=1,\dots,n-1~ \big > \longrightarrow {\mathbb C}[SL(n)]^{N^+}
\longrightarrow {\mathbb C}[N^-] \longrightarrow 0.
$$
($<w_k>$ denotes the ideal generated by the $w_k$).

Or :
$$
0 \longrightarrow \big< ~~\begin{array}{|c|} \hline
1 \\
\hline
2 \\
\hline
\vdots \\
\hline
k \\
\hline
\end{array}-1~~ \big >~+~{\mathcal {PL}}~=~{\mathcal{ PL}}_{red} \longrightarrow Sym^\bullet(\bigwedge\C^n)
\longrightarrow {\mathbb C}[N^-] \longrightarrow 0.$$

For instance, in $SL(3)$, the Pl\"ucker relation becomes in ${\mathcal {PL}}_{red}$ :
$$
\renewcommand{\arraystretch}{0.7}{
\framebox{3}} - \renewcommand{\arraystretch}{0.7}{\begin{array}{l}
\framebox{2}\\
\framebox{3} \\
\end{array}}+\renewcommand{\arraystretch}{0.7}{\begin{array}{l}
\framebox{1}\framebox{2}\\
\framebox{3} \\
\end{array}}=0.$$

Now, we look for a basis for $\C[N^-]$, by selecting some semi standard Young tableaux.\\

\begin{defn}

\n

Let $T$ be a semi standard Young tableau such that its first column
begins by $\begin{array}{|c|} \hline
1 \\
\hline
2 \\
\hline
\vdots \\
\hline
k \\
\hline
\vdots\\\hline\end{array}.$
\n We say that we can "push" $T$ if we shift the $k$ firsts rows of
$T$ to the left and we delete the column $\begin{array}{|c|} \hline
1 \\
\hline
2 \\
\hline
\vdots \\
\hline
k \\
\hline
\end{array}$ which spill out. We denote by $P(T)$ the new tableau
obtained. A tableau $T$ is said {\bf{quasi standard}} if $T$ is a semi standard Young tableau and $P(T)$
is not a semi standard tableau.
\end{defn}

\begin{exple} The $\mathfrak{sl}(3)$ case

\n

The tableaux

\n

$$\renewcommand{\arraystretch}{0.7}{\begin{array}{l}
\framebox{2}\framebox{1}\\
\framebox{3} \\
\end{array}},~~~\renewcommand{\arraystretch}{0.7}{\begin{array}{l}
\framebox{1}\framebox{3}\\
\framebox{2} \\
\end{array}} ~~\hbox{and} ~~\renewcommand{\arraystretch}{0.7}{\begin{array}{l}
\framebox{1}\framebox{2}\\
\framebox{3} \\
\end{array}}$$
are not quasi standards, but the tableaux
$$\framebox{3}~~~~\hbox{ and}~~\renewcommand{\arraystretch}{0.7}{\begin{array}{l}
\framebox{2}\\
\framebox{3} \\
\end{array}}$$
are quasi standard.\\
\end{exple}

To find a basis of ${\mathbb C}[N^-]$, adapted to representations of
$SL(n)$, we restrict ourselves to quasi standard Young tableaux.

\begin{thm}

\n

The set of quasi standard Young tableaux form a basis for the
algebra ${\mathbb C}[N^-].$
\end{thm}

Let us denote $\pi$ the canonical mapping :
$$
\pi~:~\mathbb S^\bullet=Sym^\bullet(\bigwedge\C^n)/\mathcal{PL}~\longrightarrow~
\mathbb C[N^-]=Sym^\bullet(\bigwedge\C^n)/\mathcal{PL}_{red}.
$$

The algebra of polynomial functions on $N^-$ is an indecomposable
$N^+$-module. Indeed, the action of $N^+$ on $\mathbb C[N^-]$ is defined by :
$$
n^+f(n^-)=f(~^tn^+n^-).
$$

Each module $\pi({\mathbb S}^{\lambda})={\mathbb S}^{\lambda}_{~| N^{-}}$ is
occurring in ${\mathbb C}[N^-]$.

$\pi({\mathbb S}^{\lambda})$ is generated by the lowest weight vector $\pi(w_\lambda)$ (we say it is a
monogenic module). Suppose $\lambda=a_1\omega_1+\dots+a_{n-1}\omega_{n-1}$, for each simple root $\alpha_i$, we
have :
$$
X_{\alpha_i}^{a_i}\pi(w_\lambda)\neq0~~~~\hbox{ and }~~X_{\alpha_i}^{a_i+1}\pi(w_\lambda)=0.
$$
Moreover, if $W$ is any monogenic, finite dimensional, $N^+$ module generated
by a vector $w$, for which the $X_{\alpha_i}$ are nilpotent, and if the integral numbers $a_i$ are defined by :
$$
X_{\alpha_i}^{a_i}w\neq0~~~~\hbox{ and }~~X_{\alpha_i}^{a_i+1}w=0,
$$
then $W$ is isomorphic to a quotient of $\pi({\mathbb S}^{\lambda})$.
The result in \cite{[ABW]} is :

\begin{prop}

\n

A parametrization of a basis for the quotient $\pi(\mathbb S^\lambda)={\mathbb
S}^{\lambda}_{~| N^{+}}$ is given by the set of quasi standard Young
tableaux of shape $\leq\lambda.$
\end{prop}

\begin{exple}

\n

For the Lie algebra $\mathfrak{sl}(3),$ we get the picture :
\vskip3cm
\begin{center}
\begin{picture}(180,50)(230,200) \toutpetit

\path(175,150)(335,150) \path(335,150)(415,30) \path(255,30)(415,30)
\path(175,150)(255,30) \path(215,90)(375,90) \path(335,30)(255,150)
\path(295,90)(255,150) \path(375,90)(335,30) \path(295,90)(335,30)
\path(295,90)(335,150) \path(295,90)(335,150) \path(255,150)(295,90)
\path(295,90)(255,30) \path(255,150)(215,90) \path(145,200)(175,150)
\path(410,30)(450,30)

\put(175,150){\circle{3}} \put(335,150){\circle{3}}
\put(415,30){\circle{3}} \put(255,30){\circle{3}}
\put(215,90){\circle{3}} \put(375,90){\circle{3}}
\put(335,30){\circle{3}} \put(295,90){\circle{3}}
\put(255,150){\circle{3}} \put(375,90){\circle{3}}

\put(198,77) {$\begin{array}{l}
\framebox{1}\\
\framebox{3}\\
\end{array}$}

\put(282,72){$\begin{array}{l}
\framebox{2}\\
\framebox{3}\\
\end{array}$}

\put(241,163){$\begin{array}{l}
\framebox{1}\framebox{2}\\
\framebox{3}\framebox{3}\\
\end{array}$}

\put(141,148){$\begin{array}{l}
\framebox{1}\framebox{1}\\
\framebox{3}\framebox{3}\\
\end{array}$}

\put(337,150){$\begin{array}{l}
\framebox{2}\framebox{2}\\
\framebox{3}\framebox{3}\\
\end{array}$}

\put(313.5,163){$\begin{array}{l}
\framebox{2}\framebox{3}\\
\framebox{3}\\
\end{array}$}

\put(314,136){$\begin{array}{l}
\framebox{3}\framebox{3}\\
\end{array}$}

\put(328,19.5){\framebox{2}}

\put(292,101){\framebox{3}}

\put(403,19.5){$\begin{array}{l}
\framebox{2}\framebox{2}\\
\end{array}$}

\put(374,95){$\begin{array}{l}
\framebox{2}\framebox{2}\\
\framebox{3}\\
\end{array}$}

\put(345,80){$\begin{array}{l}
\framebox{2}\framebox{3}\\
\end{array}$}

\put(255,22){0}

\end{picture}
\end{center}
\end{exple}

\vskip 7cm

\section{Principle of our construction. Fundamental representations}

\vskip0.7cm

The purpose of this article is to adress in the same way the rank
two semisimple Lie algebras. We will describe first the semi standard Young tableaux for the algebras $A_1 \times A_1, A_2,
C_2$ and $G_2$ then the quasi standard Young tableaux for these algebras. For this, we start to realize the rank
two semisimple Lie algebras as subalgebras of $\mathfrak{sl}(n)$ for $n=4,3, 4, 7$ in such a way that the simples
coroots $X_{-\alpha}$ and $X_{-\beta} (\alpha$ denotes the first "short" simple root and $\beta$ denotes the second "long"
simple root$)$ are matrices such as :
$$\left.\aligned t X_{-\alpha}& \longmapsto ~~\hbox{first ~~column ~~of}~~ t X_{-\alpha}\\
t X_{-\alpha}+ s X_{-\beta}& \longmapsto ~~\hbox{two firsts ~~columns ~~of}~~ t X_{-\alpha}+ s X_{-\beta}\endaligned\right.\eqno{(*)}$$
are one-to-one.\\
Explicitely, we take the following realizations:\\

\n

$\hskip-0.7cm \underline{\underline{A_1 \times A_1= \mathfrak{sl}(2) \times \mathfrak{sl}(2)}}:$\\
Let $(g_1,g_2) \in \mathfrak{sl}(2) \times \mathfrak{sl}(2)$ where
$g_i= \left(\begin{matrix}a_i&b_i\\c_i&d_i\end{matrix}\right)$ such that $a_i+ d_i=0$. We thus modify the
natural realization of the Lie algebra $A_1 \times A_1$ as :
$$X= \left(\begin{matrix}a_1&b_1&0&0\\c_1&d_1&0&0\\0&0&a_2&b_2\\0&0&c_2&d_2\end{matrix}
\right) \longmapsto \left(\begin{matrix}a_1&0&b_1&0\\0&a_2&0&b_2\\c_1&0&d_1&0\\0&c_2&0&d_2\end{matrix}
\right),$$
(we exchange basis vectors $2$ and $3$). Then
$$N^-= \left\{\left(
\begin{array}{cccc}
1 & 0 & 0 &0\\
0 & 1  & 0 & 0\\
x &0 & 1& 0\\
0& y & 0 & 1\\
 \end{array}
\right),~ x, y \in \mathbb{C}\right\}.$$
\hskip-0.4cm $\underline{\underline{A_2= \mathfrak{sl}(3)} }:$\\

\n

\n Let $ g \in \mathfrak{sl}(3)$ i.e $$g= \left(\begin{matrix}a_1&b_1&c_1\\a_2&b_2&c_2\\a_3&b_3&c_3\end{matrix}
\right)~~ \hbox{such ~~that}~~ a_1+ b_2+ c_3=0.$$
then $$N^-= \left\{\left(
\begin{array}{ccc}
1 & 0 & 0\\
x & 1  & 0\\
z &y & 1\\
 \end{array}
\right),~ x, y, z \in \mathbb{C}\right\}.$$

With this parametrization, we immediately see the Pl\"ucker relation in ${\mathcal{ PL}}_{red}$ :
$$
\renewcommand{\arraystretch}{0.7}{
\framebox{3}}~(g) - \renewcommand{\arraystretch}{0.7}{\begin{array}{l}
\framebox{2}\\
\framebox{3} \\
\end{array}}~(g) +\renewcommand{\arraystretch}{0.7}{\begin{array}{l}
\framebox{1}\framebox{2}\\
\framebox{3} \\
\end{array}}~(g)=z-(xy-z)+yx=0.$$

\hskip-0.4cm $\underline{\underline{C_2= \mathfrak{sp}(4)} }:$\\

\n

\hskip-0.4cm The natural realization of the Lie algebra $\mathfrak{sp}(4)$ is given by $X=\left(\begin{matrix}A&B\\ C&-~^tA\end{matrix}\right)$
with $A$, $B$, $C$ $2\times2$ matrices, and $~^tB=B$ and $~^tC=C$. We modify this realization by permuting the
basis vectors 3 and 4 :
$$
X= \left(\begin{matrix}a&b&u&v\\c&d&v&w\\x&y&-a&-c\\y&z&-b&-d\end{matrix}
\right) \longmapsto \left(\begin{matrix}a&b&v&u\\c&d&w&v\\y&z&-d&-b\\x&y&-c&-a\end{matrix}
\right).
$$

\n

Then the group $N^-$ becomes :
$$
N^-= \left\{\left(
\begin{array}{cccc}
1 & 0 & 0 &0\\
x & 1  & 0 & 0\\
z &u & 1& 0\\
y& z-xu & -x & 1\\
 \end{array}
\right),\quad x,~y,~z,~u\in\C\right\}.
$$

\hskip-0.4cm $\underline{\underline{G_2}}:$\\

\n

\hskip-0.4cm The natural realization of the Lie algebra $G_2$ is given by :
$$
X= \left(\begin{matrix}A&V&-j(\frac{W}{\sqrt{2}})\\^{-t}W&0&^{-t}V\\-j(\frac{V}{\sqrt{2}})&W&^{-t}A
\end{matrix}\right)
$$
where $V$, $W$ are $3\times1$ column-matrices, $j(U)$ is the $3\times3$ matrix of the exterior product in $\C^3$ :
$j(U)V=U \wedge V$ and $A$ is a $3\times3$ matrix such that $tr(A)=0$.

\n To imbed $\mathfrak{n}^-$ in the space of lower triangular
matrices, we effect the permutation
$\left(\begin{matrix}1&2&3&4&5&6&7\\7&2&1&4&5&6&3
\end{matrix}\right)$ on the vector basis. Then, we obtain the Lie algebra :
$$\mathfrak{n}^-=\left\{
\left(\begin{matrix}0&0&0&0&0&0&0\\-x&0&0&0&0&0&0\\y&a&0&0&0&0&0\\\sqrt{2}z&\sqrt{2}y&\sqrt{2}x&0&0&0&0\\
-b&-z&0&-\sqrt{2}x&0&0&0\\ -c&0&z& -\sqrt{2}y&-a&0&0\\
0&c&b&-\sqrt{2}z&-y&x&0
\end{matrix}
\right)\right\}$$
and the following corresponding group : $N^-$ is the set of matrices :
$$
\left(\begin{matrix}1&0&0&0&0&0&0\\x&1&0&0&0&0&0\\y&a&1&0&0&0&0\\z&-\sqrt{2}ax+\sqrt{2}y&-\sqrt{2}x&1&0&0&0\\
b&-ax^2+xy-\frac{\sqrt{2}}{2}z&-x^2&\sqrt{2}x&1&0&0\\ c&axy+\frac{\sqrt{2}}{2}az-y^2&xy+\frac{\sqrt{2}}{2}z& -\sqrt{2}y&-a&1&0\\
-yb-xc-\frac{z^2}{2}&\frac{\sqrt{2}}{2}axz-ab-\frac{\sqrt{2}}{2}yz-c&\frac{\sqrt{2}}{2}xz-b&-z&-y+ax&-x&1
\end{matrix}
\right),$$
with $a$, $b$, $c$, $x$, $y$, $z$ in $\C$.\\

In each case, we consider first the Young tableaux with 1 column and 1 or 2 rows, corresponding to particular
subrepresentations in $\C^n$ ($n=4,3,4,7$) and $\wedge^2\C^n$, which are isomorphic to the fundamental
representations $\Gamma_{1,0}$ and $\Gamma_{0,1}$ of the Lie algebra. This selection of tableaux can be
viewed as the traduction of some ``internal'' Pl\"ucker relations for our Lie algebra.\\

\hskip-0.4cm $\underline{\underline{A_1\times A_1= \mathfrak{sl}(2)\times\mathfrak{sl}(2)}}:$\\

\noindent
The $\Gamma_{1,0}$ representation occurs in $\C^4$, we find the basis $\framebox{1}~$, $\framebox{3}$ and 2
internal Pl\"ucker relations
$$
\framebox{2}=0,\quad\framebox{4}=0.
$$

The $\Gamma_{0,1}$ representation occurs in $\wedge^2\C^4$, we find the basis $\renewcommand{\arraystretch}{0.7}
\begin{array}{l}\framebox{1}\\
\framebox{2}\end{array}$, and $\renewcommand{\arraystretch}{0.7}
\begin{array}{l}\framebox{1}\\
\framebox{4}\end{array}$
and 4 internal Pl\"ucker relations
$$
\renewcommand{\arraystretch}{0.7}{\begin{array}{l}\framebox{2}\\
\framebox{4}\end{array}=0},\quad\renewcommand{\arraystretch}{0.7}{\begin{array}{l}\framebox{1}\\
\framebox{3}\end{array}=0},\quad\renewcommand{\arraystretch}{0.7}{\begin{array}{l}\framebox{2}\\
\framebox{3}\end{array}}=-\framebox{3}\quad\text{ and }~~\renewcommand{\arraystretch}{0.7}{\begin{array}{l}\framebox{3}\\
\framebox{4}\end{array}=-\begin{array}{l}\framebox{1}\framebox{3}\\
\framebox{4}\end{array}}.
$$
Thus we get the following Young semi standard tableaux with 1 column, for $\mathfrak{sl}(2)\times\mathfrak{sl}(2)$ :
$$
\framebox{1}~,~~\framebox{3}~,~~\renewcommand{\arraystretch}{0.7}{\begin{array}{l}\framebox{1}\\
\framebox{2}\end{array}},~\text{ and }~\renewcommand{\arraystretch}{0.7}{\begin{array}{l}\framebox{1}\\
\framebox{4}\end{array}}.
$$

\hskip-0.4cm $\underline{\underline{A_2= \mathfrak{sl}(3)}}:$\\

\n
By definition, there is no internal Pl\"ucker relations for $A_2$, the semi standard Young tableaux with 1 column are :
$$
\framebox{1}~,~~\framebox{2}~,~~\framebox{3}~,~~\renewcommand{\arraystretch}{0.7}{\begin{array}{l}\framebox{1}\\
\framebox{2}\end{array}},~~\renewcommand{\arraystretch}{0.7}{\begin{array}{l}\framebox{1}\\
\framebox{3}\end{array}},~\text{ and }~\renewcommand{\arraystretch}{0.7}{\begin{array}{l}\framebox{2}\\
\framebox{3}\end{array}}.
$$

\hskip-0.4cm $\underline{\underline{C_2= \mathfrak{sp}(4)} }:$\\

\n

The $\Gamma_{1,0}$ representation occurs in $\C^4$, we find the basis $\framebox{1}~$, $\framebox{2}~$,
$\framebox{3}~$ and $\framebox{4}$.\\

The $\Gamma_{0,1}$ representation is the quotient of $\wedge^2 \C^4$ by the invariant symplectic form.
Then we have 1 internal Pl\"ucker relation which is written as follows:
$$
\renewcommand{\arraystretch}{0.7}{\begin{array}{l}\framebox{1}\\
\framebox{4}\end{array}- \begin{array}{l}\framebox{2}\\
\framebox{3}\end{array}=0}.
$$
Thus we choose the Young semi standard tableaux with 1 column, for $\mathfrak{sp}(4)$ :
$$
\framebox{1}~,~~\framebox{2}~,~~\framebox{3}~,~~\framebox{4}~,~~\renewcommand{\arraystretch}{0.7}{\begin{array}{l}\framebox{1}\\
\framebox{2}\end{array},~~\begin{array}{l}\framebox{1}\\
\framebox{3}\end{array} ,~~\begin{array}{l}\framebox{2}\\
\framebox{3}\end{array},~~\begin{array}{l}\framebox{2}\\
\framebox{4}\end{array},~\text{ and }~\begin{array}{l}\framebox{3}\\
\framebox{4}\end{array}}.
$$
This choice does coincide with the choice made in \cite{[ADLMPPrW]}.

\hskip-0.4cm $\underline{\underline{G_2}}:$\\

\n

The $\Gamma_{1,0}$ representation occurs in $\C^7$, we find the basis $\framebox{1}~$, $\framebox{2}~$,
$\framebox{3}~$, $\framebox{4}~$, $\framebox{5}~$, $\framebox{6}~$ and $\framebox{7}$.\\

The $\Gamma_{0,1}$ representation is the quotient of $\wedge^2 \C^7$ by a seven dimensional module.
Then we have 7 internal Pl\"ucker relation which is written as follows:
$$
\renewcommand{\arraystretch}{0.7}{\begin{array}{l}\framebox{1}\\
\framebox{4}\end{array}+ \sqrt{2} \begin{array}{l}\framebox{2}\\
\framebox{3}\end{array}=0~,~~
\begin{array}{l}\framebox{2}\\
\framebox{4}\end{array}- \sqrt{2} \begin{array}{l}\framebox{1}\\
\framebox{5}\end{array}=0~,~~
\begin{array}{l}\framebox{3}\\
\framebox{4}\end{array}+ \sqrt{2} \begin{array}{l}\framebox{1}\\
\framebox{6}\end{array}=0~,~~
\begin{array}{l}\framebox{4}\\
\framebox{5}\end{array}+ \sqrt{2} \begin{array}{l}\framebox{2}\\
\framebox{7}\end{array}=0},$$
$$\renewcommand{\arraystretch}{0.7}{\begin{array}{l}\framebox{4}\\
\framebox{6}\end{array}- \sqrt{2} \begin{array}{l}\framebox{3}\\
\framebox{7}\end{array}=0~,~~
\begin{array}{l}\framebox{4}\\
\framebox{7}\end{array}+ \sqrt{2} \begin{array}{l}\framebox{5}\\
\framebox{6}\end{array}=0}~,~\hbox{and}
\renewcommand{\arraystretch}{0.7}{\begin{array}{l}\framebox{1}\\
\framebox{7}\end{array}- \begin{array}{l}\framebox{2}\\
\framebox{6}\end{array}- \begin{array}{l}\framebox{3}\\
\framebox{5}\end{array}=0}.
$$
Indeed, in view of the lower triangular matrices in $G_2$, with 1 on the diagonal, we find directly those relations are holding
for the corresponding functions. Moreover, these relations are covariant under the action of the diagonal matrices,
there are holding for the corresponding functions on the lower triangular matrices in $G_2$, with any non vanishing
diagonal entries, thus by $N^+$ invariance, they hold on $G_2$.\\

Thus we choose the Young semi standard tableaux with 1 column, for $G_2$ :
$$
\framebox{1}~,~~\framebox{2}~,~~\framebox{3}~,~~\framebox{4}~,~~\framebox{5}~,~~\framebox{6}~,~~\framebox{7}~$$
$$\renewcommand{\arraystretch}{0.7}{\begin{array}{l}\framebox{1}\\
\framebox{2}\end{array}, \begin{array}{l}\framebox{1}\\
\framebox{3}\end{array}, \begin{array}{l}\framebox{1}\\
\framebox{4}\end{array}, \begin{array}{l}\framebox{1}\\
\framebox{5}\end{array}, \begin{array}{l}\framebox{1}\\
\framebox{6}\end{array}, \begin{array}{l}\framebox{1}\\
\framebox{7}\end{array}, \begin{array}{l}\framebox{2}\\
\framebox{5}\end{array}, \begin{array}{l}\framebox{2}\\
\framebox{6}\end{array}, \begin{array}{l}\framebox{2}\\
\framebox{7}\end{array}, \begin{array}{l}\framebox{3}\\
\framebox{6}\end{array}, \begin{array}{l}\framebox{3}\\
\framebox{7}\end{array}, \begin{array}{l}\framebox{4}\\
\framebox{7}\end{array}, \begin{array}{l}\framebox{5}\\
\framebox{7}\end{array}, \text{ and } \begin{array}{l}\framebox{6}\\
\framebox{7}\end{array}}.
$$
This choice does coincide with the choice made in \cite{[ADLMPPrW]}.
\section{Semi standard Young tableaux for the rank two semisimple Lie algebras}

\vskip0.7cm

Following \cite{[ADLMPPrW]}, we have a construction of semi standard Young tableaux for $\Gamma_{a,b}$,
for any $a$ and $b$, knowing those of $\Gamma_{0,1}$ and $\Gamma_{1,0}$. In fact, by a general result
of Kostant (see \cite{[FH]} for instance), each non semi standard Young tableau contains a non semi standard tableau with 2 columns.
Thus, it is sufficient to determine all non semi standard tableaux with 2 columns. (In fact we
shall get conditions of $1$ or $2$ succesive columns $T^{(i)}$ and $ T^{(i+1)}$ in the tableau $T$).\\

We begin to look the fundamental representations $\Gamma_{0,1}$ and $\Gamma_{1,0}$ for the rank two semisimple
Lie algebras as spaces generated by a succession of action of $X_{-\alpha}$ and $X_{-\beta}$ on the highest weight
vector.\\

\hskip-0.4cm $\underline{\underline{A_1\times A_1}}:$

\vskip0.3cm

The fundamental representations look like :
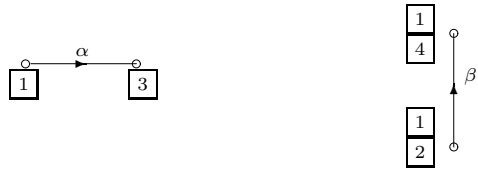
\begin{center}
\begin{picture}(50,20)(200,130)\toutpetit

\path(150,120)(190,120)

\put(170,120){\vector(1,0){1}}

\put(167,123){$\alpha$}

\put(190,120){\circle{3}}

\put(148,120){\circle{3}}

\put(142,110) {\framebox{1}}

\put(187,110){\framebox{3}}
\end{picture}
\end{center}
\begin{center}
\begin{picture}(40,70)(0,0) \toutpetit

\path(105,32)(105,75)

\put(105,55){\vector(0,1){1}}

\put(109,57){$\beta$}

\put(105,32){\circle{3}}

\put(105,75){\circle{3}}

\put(82,34){$\begin{array}{l}
 \framebox{1}  \\
  \framebox{2}  \\
\end{array}$}

\put(82,73){$\begin{array}{l}
  \framebox{1}    \\
  \framebox{4}   \\
\end{array}$}
\end{picture}
\end{center}

We associate to these drawing the two following ordered sets (respectively) :
\begin{center}
\begin{picture}(50,10)(250,150)\toutpetit
\put(190,128){\circle{3}}

\put(188,121){$\alpha$}
\end{picture}
\end{center}
\begin{center}
\begin{picture}(50,0)(170,57) \toutpetit

\put(195,50){\circle{3}}

\put(193,41){$\beta$}
\end{picture}
\end{center}
\vskip0.8cm
$\underline{\underline{A_2}}:$\\

The fundamental representations look like :
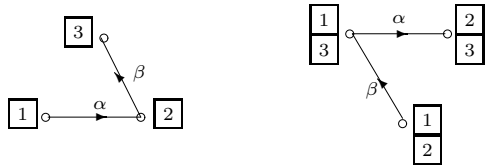
\begin{center}
\begin{picture}(180,50)(190,130) \toutpetit

\path(195,120)(230,120)

\put(216,120){\vector(1,0){1}}

\put(212,123){$\alpha$}

\put(230,120){\circle{3}}

\put(194,120){\circle{3}}

\put(180,119) {\framebox{1}}

\put(235,119){\framebox{2}}

\path(230,120)(215,150)

\put(221,136){\vector(-1,1){1}}

\put(227,135){$\beta$}

\put(216,150){\circle{3}}

\put(201,150){\framebox{3}}
\end{picture}
\end{center}

\begin{center}
\begin{picture}(200,100)(40,-5) \toutpetit

\path(170,120)(205,120)

\put(190,120){\vector(1,0){1}}

\put(185,124){$\alpha$}

\put(206,120){\circle{3}}

\put(169,120){\circle{3}}

\put(149,118) {$\begin{array}{l}
\framebox{1}\\
\framebox{3} \\
\end{array}$}

\put(204,118){$\begin{array}{l}
\framebox{2}\\
\framebox{3} \\
\end{array}$}

\path(170,120)(189,88)

\put(189,86){\circle{3}}

\put(180.5,102){\vector(-1,1){1}}

\put(175,97){$\beta$}

\put(188,80){$\begin{array}{l}
\framebox{1}\\
\framebox{2} \\
\end{array}$}
\end{picture}
\end{center}
\vskip-2cm Then, we represent these drawing by the two following ordered sets (respectively):
\begin{center}
\begin{picture}(180,90)(70,30)\toutpetit

\path(105,35)(105,100)

\put(105,82){\vector(0,-1){1}}

\put(105,52){\vector(0,-1){1}}

\put(108,83){$\alpha$}
\put(108,51){$\beta$}

\put(105,65){\circle{3}}
\put(105,100){\circle{3}}
\put(105,34){\circle{3}}

\put(90,32){\framebox{3}}

\put(90,97){\framebox{1}}

\put(90,65){\framebox{2}}
\end{picture}
\end{center}
\begin{center}
\begin{picture}(10,-230)(10,15)\toutpetit

\path(105,35)(105,100)

\put(104.5,82){\vector(0,-1){1}}

\put(104.5,52){\vector(0,-1){1}}

\put(109,83){$\beta$}

\put(109,51){$\alpha$}

\put(104.5,65){\circle{3}}

\put(104.5,100){\circle{3}}

\put(104.5,34){\circle{3}}

\put(84,32){$\begin{array}{l}
\framebox{2}\\
\framebox{3} \\
\end{array}$}

\put(84,97){$\begin{array}{l}
\framebox{1}\\
\framebox{2} \\
\end{array}$}

\put(84,65){$\begin{array}{l}
\framebox{1}\\
\framebox{3} \\
\end{array}$}
\end{picture}
\end{center}

\vfill\eject
 $\underline{\underline{C_2}}:$\\

The fundamental representations look like :

\begin{center}
\begin{picture}(90,50)(250,130)\toutpetit

\path(176,128)(200,106)

\path(224,128)(200,153)

\path(200,106)(200,153)

\put(224,128){\circle{3}}

\put(176,128){\circle{3}}

\put(162,125){\framebox{1}}

\put(228,125){\framebox{4}}

\put(186,157){\framebox{3}}

\put(194,95){\framebox{2}}

\put(200,106){\circle{3}}

\put(217,136){\vector(1,-1){1}}

\put(200,130){\vector(0,1){1}}

\put(219,140){$\alpha$}

\put(200,153){\circle{3}}

\put(186,119){\vector(1,-1){1}}

\put(178,115){$\alpha$}

\put(205,132){$\beta$}
\end{picture}
\end{center}
\begin{center}
\begin{picture}(160,70)(0,0)\toutpetit
\path(105,54)(105,99)

\put(105,78){\vector(0,1){1}}

\put(98,76){$\beta$}

\put(105,100){\circle{3}}

\put(105,54){\circle{3}}

\path(148,54)(148,99)

\put(147,80){\vector(0,1){1}}

\put(152,78){$\beta$}

\put(148,100){\circle{3}}

\path(149,54)(105,100)

\put(149,54){\circle{3}}

\put(127,77){\circle{3}}

\put(138,67){\vector(1,-1){1}}

\put(118,90){$\alpha$}

\put(116,89){\vector(1,-1){1}}

\put(139,70){$\alpha$}

\put(86,103){$\begin{array}{l}
 \framebox{1}  \\
  \framebox{3}  \\
\end{array}$}

\put(147,47){$\begin{array}{l}
 \framebox{2}  \\
  \framebox{4}  \\
\end{array}$}

\put(147,103){$\begin{array}{l}
\framebox{3}  \\
\framebox{4}  \\
\end{array}$}

\put(109,65){$\begin{array}{l}
 \framebox{2}  \\
  \framebox{3}  \\
\end{array}$}

\put(86,47){$\begin{array}{l}
 \framebox{1}  \\
  \framebox{2}  \\
\end{array}$}
\end{picture}
\end{center}
\vskip-0.5cm Then, We associate to these drawing the two following ordered sets (respectively) :
\begin{center}
\begin{picture}(200,52)(100,90)\toutpetit

\path(105,0)(105,100)

\put(105,75){\vector(0,-1){1}}

\put(105,43){\vector(0,-1){1}}

\put(108,77){$\alpha$}

\put(108,45){$\beta$}

\put(105,30){\circle{3}}

\put(105,99){\circle{3}}

\put(105,62){\circle{3}}

\put(89,30) {\framebox{3}}

\put(89,98) {\framebox{1}}

\put(89,62) {\framebox{2}}

\put(105,1){\circle{3}}

\put(89,-1) {\framebox{4}}

\put(108,17){$\alpha$}

\put(105,16){\vector(0,-1){1}}
\end{picture}
\end{center}
\begin{center}
\begin{picture}(100,90)(0,30)\toutpetit

\path(105,5)(105,140)

\put(105,88){\vector(0,-1){1}}

\put(105,54){\vector(0,-1){1}}

\put(109,90){$\alpha$}

\put(109,55){$\alpha$}

\put(105,5){\circle{3}}

\put(105,141){\circle{3}}

\put(105,108){\circle{3}}

\put(83,37) {$\begin{array}{l}
\framebox{2}\\
\framebox{4} \\
\end{array}$}

\put(83,107) {$\begin{array}{l}
\framebox{1}\\
\framebox{3} \\
\end{array}$}

\put(83,74) {$\begin{array}{l}
\framebox{2}\\
\framebox{3} \\
\end{array}$}

\put(105,74){\circle{3}}

\put(83,5) {$\begin{array}{l}
\framebox{3}\\
\framebox{4} \\
\end{array}$}

\put(109,23){$\beta$}

\put(105,23){\vector(0,-1){1}}

\put(105,38){\circle{3}}

\put(83,139) {$\begin{array}{l}
\framebox{1}\\
\framebox{2} \\
\end{array}$}

\put(109,125){$\beta$}

\put(105,123){\vector(0,-1){1}}
\end{picture}
\end{center}
\vskip1cm
 $\underline{\underline{G_2}}:$\\

The fundamental representations look like :
\newpage
\begin{center}
\begin{picture}(52,120)(350,90)\toutpetit
\path(190,80)(237,80)

\path(237,158)(189,158)

\put(168,119){\circle{3}}

\put(189,80){\circle{3}}

\put(237,80){\circle{3}}

\put(260,119){\circle{3}}

\put(188,159){\circle{3}}

\put(237,159){\circle{3}}

\put(219,80){\vector(1,0){1}}

\put(217,158){\vector(1,0){1}}

\put(217,74){$\alpha$}

\put(215,162){$\alpha$}

\path(168,119)(258,119)

\put(212,119){\circle{3}}

\path(260,119)(188,158)

\path(168,118)(237,80)

\put(208,96){\vector(-1,1){1}}

\put(195,119){\vector(1,0){1}}

\put(234,119){\vector(1,0){1}}

\put(205,149){\vector(-1,1){1}}

\put(203,92){$\beta$}

\put(190,122){$\alpha$}

\put(229,122){$\alpha$}

\put(203,141){$\beta$}

\put(207,110) {\framebox{4}}

\put(241,159) {\framebox{7}}

\put(175,159){\framebox{6}}

\put(155,117){\framebox{3}}

\put(175,77){\framebox{1}}

\put(262,117){\framebox{5}}

\put(241,77){\framebox{2}}

\end{picture}
\end{center}
\begin{center}
\begin{picture}(52,0)(140,150) \toutpetit

\path(378,210)(275,265)

\path(275,33)(172,88)

\put(172,88){\circle{3}}

\put(172,210){\circle{3}}

\put(275,265){\circle{3}}

\put(378,210){\circle{3}}

\put(378,88){\circle{3}}

\put(275,33){\circle{3}}

\path(241,210)(380,210)

\path(172,88)(378,88)

\put(241,88){\circle{3}}

\put(310,88){\circle{3}}

\put(241,210){\circle{3}}

\put(310,210){\circle{3}}

\path(274,149)(345,149)

\put(206,149){\circle{3}}

\put(343,149){\circle{3}}

\put(274,149){\circle{3}}

\path(378,88)(274,149)

\path(344,149)(241,210)

\bezier{20}(206,149)(237,180)(274,149)

\put(238,165){\vector(1,0){1}}

\put(235,159){$\alpha$}

\bezier{20}(274,149)(305,180)(344,149)

\put(305,165){\vector(1,0){1}}

\put(302,159){$\alpha$}

\bezier{70}(206,149)(265,115)(310,88)

\put(241,128){\vector(-1,1){1}}

\put(238,121){$\beta$}

\bezier{70}(206,149)(240,149)(274,149)

\bezier{70}(274,149)(225,180)(172,210)

\put(215,186){\vector(-1,1){1}}

\put(210,182){$\alpha$}

\bezier{70}(172,210)(205,210)(241,210)

\put(205,210){\vector(1,0){1}}

\put(204,213){$\alpha$}

\put(205,88){\vector(1,0){1}}

\put(203,82){$\alpha$}

\put(276,88){\vector(1,0){1}}

\put(274,82){$\alpha$}

\put(345,88){\vector(1,0){1}}

\put(343,82){$\alpha$}

\put(274,210){\vector(1,0){1}}

\put(271,213){$\alpha$}

\put(345,210){\vector(1,0){1}}

\put(342,213){$\alpha$}

\put(325,119){\vector(-1,1){1}}

\put(327,122){$\beta$}

\put(239,149){\vector(1,0){1}}

\put(237,144){$\alpha$}

\put(300,149){\vector(1,0){1}}

\put(298,144){$\alpha$}

\put(284,186){\vector(-1,1){1}}

\put(288,187){$\beta$}

\put(217,65){\vector(-1,1){1}}

\put(215,57){$\beta$}

\put(319,243){\vector(-1,1){1}}

\put(323,247){$\beta$}

\put(151,89) {$\begin{array}{l}
\framebox{1}\\
\framebox{3} \\
\end{array}$}

\put(232,73) {$\begin{array}{l}
\framebox{1}\\
\framebox{4} \\
\end{array}$}

\put(300,73) {$\begin{array}{l}
\framebox{1}\\
\framebox{5} \\
\end{array}$}

\put(378,87) {$\begin{array}{l}
\framebox{2}\\
\framebox{5} \\
\end{array}$}

\put(265,18) {$\begin{array}{l}
\framebox{1}\\
\framebox{2} \\
\end{array}$}

\put(186,141) {$\begin{array}{l}
\framebox{1}\\
\framebox{6} \\
\end{array}$}

\put(262,133) {$\begin{array}{l}
\framebox{2}\\
\framebox{6} \\
\end{array}$}

\put(269,162) {$\begin{array}{l}
\framebox{1}\\
\framebox{7} \\
\end{array}$}

\put(344,141) {$\begin{array}{l}
\framebox{2}\\
\framebox{7} \\
\end{array}$}

\put(151,209) {$\begin{array}{l}
\framebox{3}\\
\framebox{6} \\
\end{array}$}

\put(231,225) {$\begin{array}{l}
\framebox{3}\\
\framebox{7} \\
\end{array}$}

\put(300,225) {$\begin{array}{l}
\framebox{4}\\
\framebox{7} \\
\end{array}$}

\put(378,209) {$\begin{array}{l}
\framebox{5}\\
\framebox{7} \\
\end{array}$}

\put(265,280) {$\begin{array}{l}
\framebox{6}\\
\framebox{7} \\
\end{array}$}
\end{picture}
\end{center}

\vskip5cm
\n Then, we associate to these drawing the two following ordered sets (respectively) :

\begin{center}
\begin{picture}(200,260)(150,-70)\toutpetit

\path(105,-100)(105,100)

\put(105,-20){\vector(0,-1){1}}

\put(108,-20){$\alpha$}

\put(105,-55){\vector(0,-1){1}}

\put(108,-56){$\beta$}

\put(105,-84){\vector(0,-1){1}}

\put(108,-83){$\alpha$}

\put(105,81){\vector(0,-1){1}}

\put(105,43){\vector(0,-1){1}}

\put(108,81){$\alpha$}

\put(108,45){$\beta$}

\put(105,31){\circle{3}}

\put(105,100){\circle{3}}

\put(105,66){\circle{3}}

\put(89,31) {\framebox{3}}

\put(89,102) {\framebox{1}}

\put(89,65) {\framebox{2}}

\put(105,-3){\circle{3}}

\put(89,-3) {\framebox{4}}

\put(108,17){$\alpha$}

\put(105,16){\vector(0,-1){1}}

\put(105,-36){\circle{3}}

\put(105,-68){\circle{3}}

\put(105,-100){\circle{3}}

\put(89,-35) {\framebox{5}}

\put(89,-68) {\framebox{6}}

\put(89,-99) {\framebox{7}}
\end{picture}
\end{center}

\n

\begin{center}
\begin{picture}(80,-600)(100,-150)\toutpetit

\path(200,87)(200,0)

\put(200,87){\circle{3}}

\put(200,0){\circle{3}}

\put(200,29){\circle{3}}

\put(200,58){\circle{3}}

\bezier{30}(200,0)(215,-15)(230,-30)

\path(170,-30)(200,0)

\put(170,-30){\circle{3}}

\put(260,-60){\circle{3}}

\put(230,-30){\circle{3}}

\path(170,-30)(230,-90)

\put(230,-90){\circle{3}}

\path(199,-120)(230,-90)

\bezier{30}(170,-90)(200,-60)(230,-30)

\bezier{30}(170,-90)(185,-105)(200,-120)

\bezier{60}(230,-30)(255,-5)(260,-60)

\bezier{60}(230,-90)(260,-115)(260,-60)

\path(200,-203)(200,-120)

\put(170,-90){\circle{3}}

\put(200,-120){\circle{3}}

\put(230,-30){\circle{3}}

\put(200,-148){\circle{3}}

\put(200,-175){\circle{3}}

\put(200,-202){\circle{3}}

\put(200,-61) {\circle{3}}

\put(200,72){\vector(0,-1){1}}

\put(205,72){$\beta$}

\put(200,41){\vector(0,-1){1}}

\put(205,41){$\alpha$}

\put(200,16){\vector(0,-1){1}}

\put(205,16){$\alpha$}

\put(215,-14){\vector(1,-1){1}}

\put(209,-19){$\beta$}

\put(180,-20){\vector(-1,-1){1}}

\put(184,-23){$\alpha$}

\put(188,-47){\vector(1,-1){1}}

\put(179,-48){$\beta$}

\put(213,-46){\vector(-1,-1){1}}

\put(217,-48){$\alpha$}

\put(255,-28){\vector(1,-1){1}}

\put(256,-26){$\alpha$}

\put(257,-85){\vector(-1,-1){1}}

\put(262,-87){$\alpha$}

\put(215,-75){\vector(1,-1){1}}

\put(210,-78){$\alpha$}

\put(185,-75){\vector(-1,-1){1}}

\put(189,-77){$\beta$}

\put(185,-105){\vector(1,-1){1}}

\put(186,-102){$\alpha$}

\put(215,-105){\vector(-1,-1){1}}

\put(210,-102){$\beta$}

\put(200,-135){\vector(0,-1){1}}

\put(205,-133){$\alpha$}

\put(200,-164){\vector(0,-1){1}}

\put(205,-162){$\alpha$}

\put(200,-190){\vector(0,-1){1}}

\put(205,-192){$\beta$}

\put(179,89){$\begin{array}{l}
\framebox{1}\\
\framebox{2} \\
\end{array}$}

\put(179,61){$\begin{array}{l}
\framebox{1}\\
\framebox{3} \\
\end{array}$}

\put(179,32){$\begin{array}{l}
\framebox{1}\\
\framebox{4} \\
\end{array}$}

\put(179,4) {$\begin{array}{l}
\framebox{1}\\
\framebox{5} \\
\end{array}$}

\put(223,-16) {$\begin{array}{l}
\framebox{1}\\
\framebox{6} \\
\end{array}$}

\put(152,-28) {$\begin{array}{l}
\framebox{2}\\
\framebox{5} \\
\end{array}$}

\put(190,-45) {$\begin{array}{l}
\framebox{2}\\
\framebox{6} \\
\end{array}$}

\put(259,-58) {$\begin{array}{l}
\framebox{1}\\
\framebox{7} \\
\end{array}$}

\put(223,-105) {$\begin{array}{l}
\framebox{2}\\
\framebox{7} \\
\end{array}$}

\put(152,-94) {$\begin{array}{l}
\framebox{3}\\
\framebox{6} \\
\end{array}$}

\put(179,-127) {$\begin{array}{l}
\framebox{3}\\
\framebox{7} \\
\end{array}$}

\put(179,-151) {$\begin{array}{l}
\framebox{4}\\
\framebox{7} \\
\end{array}$}

\put(179,-176) {$\begin{array}{l}
\framebox{5}\\
\framebox{7} \\
\end{array}$}

\put(179,-206) {$\begin{array}{l}
\framebox{6}\\
\framebox{7} \\
\end{array}$}
\end{picture}
\end{center}
\vskip4cm
\n

Then, we can realize these chosen paths as the family $L$ of
ideals of some partially ordered sets $P$ (which are called {\bf{posets}} ). An  ideal in $P$ is a subset
$I \subset P$ such that if $u \in P$ and $v \leq u$, then $v \in I$. With our choice,
we take the following fundamental posets denoted $P_{1,0}$ and $P_{0,1}$ and we associate for each of them
the correspondant distributive lattice of their ideals enoted $L_{1,0}$ and $L_{0,1}$ (respectively).  \\

\n

\hskip-0.4cm $\underline{\underline{A_1 \times A_1}}:$\\

$\underline{P_{1,0}}:$  \hskip7cm $\underline{P_{0,1}}:$   \\
\begin{picture}(50,10)(120,140)\toutpetit
\put(190,128){\circle{3}}

\put(188,121){$\alpha$}
\end{picture}

\begin{picture}(50,10)(-100,120)\toutpetit
\put(190,128){\circle{3}}

\put(188,119){$\beta$}
\end{picture}
\vskip1cm
$\underline{L_{1,0}}:$  \hskip7cm $\underline{L_{0,1}}:$   \\

\begin{picture}(40,70)(30,30) \toutpetit

\path(105,32)(105,75)

\put(103,52){$\alpha$}

\put(105,30){\circle{3}}

\put(105,76){\circle{3}}

\put(91,29){$(\varnothing)$}

\put(91,75){$(\alpha)$}
\end{picture}

\begin{picture}(40,70)(-200,-40) \toutpetit

\path(105,32)(105,75)

\put(103,52){$\beta$}

\put(105,30){\circle{3}}

\put(105,76){\circle{3}}

\put(91,29){$(\varnothing)$}

\put(91,75){$(\beta)$}
\end{picture}
\vskip-1cm

$\underline{\underline{A_2}}:$\\

$\underline{P_{1,0}}:$  \hskip7cm $\underline{P_{0,1}}:$   \\

\begin{picture}(40,70)(0,-20) \toutpetit

\path(75,32)(105,2)

\put(106,8){$\beta$}

\put(73,33){\circle{3}}

\put(105,2){$\circle{3}$}

\put(79,35){$\alpha$}
\end{picture}

\begin{picture}(40,70)(-230,-90) \toutpetit

\path(75,32)(105,2)

\put(106,8){$\alpha$}

\put(73,33){\circle{3}}

\put(105,2){\circle{3}}

\put(79,35){$\beta$}
\end{picture}

\vskip-1.5cm

$\underline{L_{1,0}}:$  \hskip7cm $\underline{L_{0,1}}:$   \\

\begin{picture}(110,50)(20,50)\toutpetit

\path(105,35)(105,100)

\put(104,83){$\alpha$}

\put(104,51){$\beta$}

\put(105,65){\circle{3}}

\put(105,100){\circle{3}}

\put(105,34){\circle{3}}

\put(7,7){$\path(105,98)(117,86)$}

\put(124,92){\circle{2}}

\put(111,106){\circle{2}}

\put(129,92){$\beta$}

\put(111,109){$\alpha$}

\put(110,65){$(\beta)$}

\put(110,33){$(\varnothing)$}

\end{picture}

\begin{picture}(50,-200)(-200,40)\toutpetit

\path(105,35)(105,100)

\put(104,83){$\beta$}

\put(104,51){$\alpha$}

\put(105,65){\circle{3}}

\put(105,100){\circle{3}}

\put(105,34){\circle{3}}

\put(7,7){$\path(105,98)(117,86)$}

\put(124,93){\circle{2}}

\put(111,106){\circle{2}}

\put(129,92){$\alpha$}

\put(110,65){$(\alpha)$}

\put(111,110){$\beta$}

\put(110,33){$(\varnothing)$}

\end{picture}

$\underline{\underline{C_2}}:$\\

$\underline{P_{1,0}}:$  \hskip7cm $\underline{P_{0,1}}:$   \\

\begin{picture}(40,30)(0,20) \toutpetit

\path(75,32)(125,-22)

\put(105,5){$\beta$}

\put(75,32){\circle{3}}

\put(101,4){\circle{3}}

\put(125,-22){\circle{3}}

\put(79,35){$\alpha$}

\put(127,-19){$\alpha$}
\end{picture}

\begin{picture}(100,-150)(-150,-50)\toutpetit

\path(170,-30)(200,-60)

\path(200,-60)(170,-90)

\path(200,-120)(170,-90)

\put(170,-30){\circle{3}}

\put(200,-60){\circle{3}}

\put(170,-90){\circle{3}}

\put(200,-120){\circle{3}}

\put(170,-26){$\beta$}

\put(203,-61){$\alpha$}

\put(162,-90){$\alpha$}

\put(199,-126){$\beta$}
\end{picture}
\vskip3cm

$\underline{L_{1,0}}:$  \hskip7cm $\underline{L_{0,1}}:$   \\

\begin{picture}(110,-135)(20,120)\toutpetit

\path(105,35)(105,100)

\path(105,35)(105,3)

\put(104,21){$\alpha$}

\put(104,83){$\alpha$}

\put(104,49){$\beta$}

\put(105,65){\circle{3}}

\put(105,3){\circle{3}}

\put(105,100){\circle{3}}

\put(105,34){\circle{3}}

\put(7,7){$\path(105,98)(129,74)$}

\put(123,94){\circle{2}}

\put(135,82){\circle{2}}

\put(111,105){\circle{2}}

\put(126,97){$\beta$}

\put(111,109){$\alpha$}

\put(138,84){$\alpha$}

\put(5,-28){$\path(105,98)(117,86)$}

\put(110,70){\circle{2}}

\put(122,58){\circle{2}}

\put(112,74){$\beta$}

\put(125,59){$\alpha$}

\put(110,33){$(\alpha)$}

\put(110,3){$(\varnothing)$}
\end{picture}

\begin{picture} (100,200)(-190,-70)\toutpetit
\path(105,35)(105,100)

\path(105,35)(105,3)

$\path(105,3)(105,-25)$

\put(105,-24){\circle{3}}

\put(104,21){$\alpha$}

\put(104,83){$\beta$}

\put(104,49){$\alpha$}

\put(104,-10){$\beta$}

\put(105,65){\circle{3}}

\put(105,3){\circle{3}}

\put(105,100){\circle{3}}

\put(105,34){\circle{3}}

\put(112,106){\circle{2}}

\put(122,96){\circle{2}}

\put(112,126){\circle{2}}

\put(122,116){\circle{2}}

\put(126,97){$\beta$}

\put(122,80){\circle{2}}

\put(112,70){\circle{2}}

\put(122,60){\circle{2}}

\put(106,109){$\alpha$}

\put(125,115){$\alpha$}

\put(106,128){$\beta$}

\put(126,58){$\beta$}

\put(116,69){$\alpha$}

\put(125,78){$\alpha$}

\put(110,-24){$(\varnothing)$}

\put(-58,156){$\path(170,-30)(180,-40)$}

\put(-58,156){\path(180,-40)(170,-50)}

\put(-58,156){\path(180,-60)(170,-50)}

\put(-58,120){\path(180,-40)(170,-50)}

\put(-58,120){\path(180,-60)(170,-50)}

\put(-58,90){\path(180,-60)(170,-50)}

\put(112,40){\circle{2}}

\put(122,30){\circle{2}}

\put(110,3){$(\beta)$}

\put(115,40){$\alpha$}

\put(127,30){$\beta$}
\end{picture}

\vskip-1cm
 $\underline{\underline{G _2}}:$\\

$\underline{P_{1,0}}:$  \hskip7cm $\underline{P_{0,1}}:$   \\

\begin{picture}(100,-150)(130,0)\toutpetit

\path(170,-30)(215,-75)

\path(193,-98)(215,-75)

\path(193,-98)(238,-143)

\put(170,-30){\circle{3}}

\put(215,-75){\circle{3}}

\put(193,-98){\circle{3}}

\put(238,-143){\circle{3}}

\put(193,-53){\circle{3}}

\put(215,-120){\circle{3}}

\put(171.7,-27){$\alpha$}

\put(218,-74){$\alpha$}

\put(233,-147.6){$\alpha$}

\put(208,-124.5){$\beta$}

\put(195,-50.6){$\beta$}

\put(187.4,-102.6){$\alpha$}
\end{picture}

\begin{picture}(70,-100)(-100,20)\toutpetit

\path(170,-30)(215,-75)

\path(192,-98)(215,-75)

\path(192,-98)(215,-120)

\path(237,-53)(215,-75)

\path(237,-53)(192,-8)

\path(192,-8)(170,-30)

\path(192,-8)(215,14)

\path(192,36)(215,14)

\path(192,-54)(215,-30)

\put(170,-30){\circle{3}}

\put(215,-75){\circle{3}}

\put(192.5,-98){\circle{3}}

\put(214.5,-119.5){\circle{3}}

\put(237,-53){\circle{3}}

\put(192,-8){\circle{3}}

\put(192.5,36){\circle{3}}

\put(214.5,14){\circle{3}}

\put(192.5,-53.5){\circle{3}}

\put(214,-30){\circle{3}}

\put(195,39.5){$\beta$}

\put(218.5,15){$\alpha$}

\put(197,-8){$\alpha$}

\put(164,-29){$\alpha$}

\put(217,-28){$\beta$}

\put(185.5,-56.7){$\beta$}

\put(240,-51){$\alpha$}

\put(219,-78){$\alpha$}

\put(196,-97){$\alpha$}

\put(218,-119){$\beta$}
\end{picture}
\vskip5cm
$\underline{L_{1,0}}:$  \hskip7cm $\underline{L_{0,1}}:$   \\

\begin{picture} (100,200)(30,10)\toutpetit
\path(105,35)(105,147)

\path(105,35)(105,3)

\path(105,3)(105,-25)

\path(105,-140)(105,-25)

\put(105,-140){\circle{3}}

\put(105,147){\circle{3}}

\put(105,-94){\circle{3}}

\put(105,-48){\circle{3}}

\put(105,21){$\alpha$}

\put(105,70){$\beta$}

\put(105,-20){$\alpha$}

\put(105,-2){\circle{3}}

\put(105,98){\circle{3}}

\put(105,47){\circle{3}}

\put(110,-93){$(\beta)$}

\put(-58,10){\path(180,-60)(170,-50)}

\put(121,-49){\circle{2}}

\put(113,-40.5){\circle{2}}

\put(110,-139){$(\varnothing)$}

\put(-58,53){\path(190,-70)(170,-50)}

\put(113,2){\circle{2}}

\put(132,-17){\circle{2}}

\put(123,-8){\circle{2}}

\put(114.5,3){$\alpha$}

\put(124,-6){$\alpha$}

\put(135,-18){$\beta$}

\put(-58,107){\path(190,-70)(170,-50)}

\put(-58,107){\path(180,-39)(170,-50)}

\put(131,38){\circle{2}}

\put(121,48){\circle{2}}

\put(113,57){\circle{2}}

\put(121,67){\circle{2}}

\put(134,38){$\alpha$}

\put(124,51){$\beta$}

\put(116,58){$\alpha$}

\put(124,70){$\alpha$}

\put(104,-118){$\alpha$}

\put(104,-75){$\beta$}

\put(112,-35){$\beta$}

\put(124,-53){$\alpha$}

\put(-58,157){\path(190,-70)(170,-50)}

\put(-58,157){\path(180,-39)(170,-50)}

\put(-58,157){\path(180,-39)(170,-29)}

\put(131.5,88){\circle{2}}

\put(112.5,107){\circle{2}}

\put(122,118){\circle{2}}

\put(112.5,127.5){\circle{2}}

\put(122.5,96.5){\circle{2}}

\put(134,87){$\alpha$}

\put(124,118){$\alpha$}

\put(115,105.7){$\alpha$}

\put(116.7,93){$\beta$}

\put(115,129){$\beta$}

\put(-58,197){\path(190,-60)(170,-40)}

\put(-58,197){\path(180,-29)(170,-40)}

\put(-58,197){\path(180,-29)(160,-9)}

\put(131.7,137.6){\circle{2}}

\put(112.6,157){\circle{2}}

\put(121.4,168.5){\circle{2}}

\put(112,178){\circle{2}}

\put(121.5,147.5){\circle{2}}

\put(102.3,187.7){\circle{2}}

\put(135,140){$\alpha$}

\put(106,158){$\alpha$}

\put(122,152){$\beta$}

\put(124,170){$\alpha$}

\put(113,183){$\beta$}

\put(103,192){$\alpha$}

\put(104,125){$\alpha$}
\end{picture}

\begin{picture}(100,200)(-100,-250)\toutpetit

\path(200,87)(200,0)

\put(200,87){\circle{3}}

\put(200,0){\circle{3}}

\put(200,29){\circle{3}}

\put(200,58){\circle{3}}

\path(230,-30)(260,-60)

\path(200,0)(230,-30)

\path(170,-30)(200,0)

\put(170,-30){\circle{3}}

\put(260,-60){\circle{3}}

\put(230,-30){\circle{3}}

\path(170,-30)(230,-90)

\put(230,-90){\circle{3}}

\path(199,-120)(260,-60)

\path(170,-90)(230,-30)

\path(170,-90)(200,-120)

\path(230,-30)(260,-60)

\path(230,-90)(260,-60)

\path(200,-203)(200,-120)

\put(170,-90){\circle{3}}

\put(200,-120){\circle{3}}

\put(230,-30){\circle{3}}

\put(200,-148){\circle{3}}

\put(200,-175){\circle{3}}

\put(200,-202){\circle{3}}

\put(200,-61){\circle{3}}

\put(199,72){$\beta$}

\put(199,41){$\alpha$}

\put(199,16){$\alpha$}

\put(212,-13){$\beta$}

\put(180,-19){$\alpha$}

\put(181,-44){$\beta$}

\put(213,-46){$\alpha$}

\put(246,-45){$\alpha$}

\put(247,-73){$\alpha$}

\put(213,-75){$\alpha$}

\put(179,-79){$\beta$}

\put(182,-105){$\alpha$}

\put(215,-105){$\beta$}

\put(199,-133){$\alpha$}

\put(199,-162){$\alpha$}

\put(199,-192){$\beta$}

\put(35,126){\path(170,-30)(181,-41)}

\put(35,126){\path(181,-41)(175,-47)}

\put(35,126){\path(181,-53)(175,-47)}

\put(35,126){\path(181,-41)(187,-35)}

\put(35,126){\path(181,-30)(175,-36)}

\put(35,126){\path(187,-35)(175,-24)}

\put(35,126){\path(170,-30)(175,-24)}

\put(35,126){\path(181,-18)(175,-24)}

\put(35,126){\path(181,-18)(175,-12)}

\put(209.6,114.8){\circle{2}}

\put(210,118){$\beta$}

\put(215.7,108){\circle{2}}

\put(219,111){$\alpha$}

\put(209.5,102){\circle{2}}

\put(207,105){$\alpha$}

\put(204.9,96.7){\circle{2}}

\put(200,98){$\alpha$}

\put(210.3,90){\circle{2}}

\put(205,87){$\beta$}

\put(216,97){\circle{2}}

\put(216,99){$\beta$}

\put(222,91){\circle{2}}

\put(225,93){$\alpha$}

\put(215.5,85){\circle{2}}

\put(217,82){$\alpha$}

\put(210.4,79.7){\circle{2}}

\put(205,78){$\alpha$}

\put(216,73){\circle{2}}

\put(219,72){$\beta$}

\put(5,95){\path(181,-53)(175,-47)}

\put(5,95){\path(181,-41)(187,-35)}

\put(5,95){\path(187,-35)(175,-24)}

\put(5,95){\path(170,-30)(175,-24)}

\put(5,95){\path(181,-18)(175,-24)}

\put(5,95){\path(181,-30)(175,-36)}

\put(5,95){\path(170,-30)(181,-41)}

\put(5,95){\path(181,-41)(175,-47)}

\put(186,77.4){\circle{2}}
\put(185,80){$\alpha$}
\put(180,71){\circle{2}}
\put(176,73){$\alpha$}
\put(174.8,65){\circle{2}}
\put(168,66){$\alpha$}
\put(186,66){\circle{2}}
\put(188,67){$\beta$}
\put(180,60){\circle{2}}
\put(176,56){$\beta$}
\put(192,59.8){\circle{2}}
\put(192,55){$\alpha$}
\put(181,49){\circle{2}}
\put(178,45){$\alpha$}
\put(186,54){\circle{2}}
\put(185,50){$\alpha$}
\put(185.7,42){\circle{2}}
\put(188,39){$\beta$}

\put(40,65){\path(181,-41)(187,-35)}

\put(40,65){\path(187,-35)(175,-24)}

\put(40,65){\path(170,-30)(175,-24)}

\put(40,65){\path(181,-30)(175,-36)}

\put(40,65){\path(170,-30)(181,-41)}

\put(40,65){\path(181,-41)(175,-47)}

\put(40,65){\path(181,-53)(175,-47)}

\put(214.8,41){\circle{2}}
\put(215,43){$\alpha$}
\put(210,35.5){\circle{2}}
\put(205,34){$\alpha$}
\put(220.8,35.5){\circle{2}}
\put(224,37){$\beta$}
\put(215.5,30){\circle{2}}
\put(211,26){$\beta$}
\put(227,30.3){\circle{2}}
\put(230,30){$\alpha$}
\put(221,24){\circle{2}}
\put(222,22){$\alpha$}
\put(215.5,18.5){\circle{2}}
\put(212,14){$\alpha$}
\put(221.3,12){\circle{2}}
\put(223,11){$\alpha$}

\put(5,45){\path(181,-30)(175,-36)}

\put(5,45){\path(170,-30)(181,-41)}

\put(5,45){\path(181,-41)(175,-47)}

\put(5,45){\path(181,-53)(175,-47)}

\put(5,45){\path(181,-41)(187,-35)}

\put(5,45){\path(181,-30)(187,-35)}

\put(-20,10){\path(181,-30)(175,-36)}

\put(-20,10){\path(175,-36)(181,-41)}

\put(-20,10){\path(181,-41)(175,-47)}

\put(-20,10){\path(181,-53)(175,-47)}

\put(-20,10){\path(181,-41)(187,-35)}

\put(-20,10){\path(181,-30)(187,-35)}

\put(175,15){\circle{2}}
\put(168,14){$\alpha$}
\put(185.6,15){\circle{2}}
\put(184,19){$\beta$}
\put(180,10){\circle{2}}
\put(176,6){$\beta$}
\put(191.3,10.5){\circle{2}}
\put(193,09){$\alpha$}
\put(185.2,4.8){\circle{2}}
\put(188,04){$\alpha$}
\put(180.8,-2){\circle{2}}
\put(174,-2){$\alpha$}
\put(185.5,-8){\circle{2}}
\put(179,-10){$\beta$}

\put(65,15){\path(170,-30)(181,-41)}

\put(65,15){\path(181,-41)(175,-47)}

\put(65,15){\path(181,-53)(175,-47)}

\put(65,15){\path(181,-41)(187,-35)}

\put(234,-14){\circle{2}}
\put(230,-11){$\alpha$}
\put(240,-19.3){\circle{2}}
\put(243,-16){$\beta$}
\put(252.4,-19.4){\circle{2}}
\put(255,-21){$\alpha$}
\put(246,-26.1){\circle{2}}
\put(249,-29){$\alpha$}
\put(240.4,-32){\circle{2}}
\put(237,-29){$\alpha$}
\put(247,-39){\circle{2}}
\put(250,-38){$\beta$}

\put(20,-35){\path(175,-36)(181,-41)}

\put(20,-35){\path(181,-41)(175,-47)}

\put(20,-35){\path(181,-53)(175,-47)}

\put(20,-35){\path(181,-41)(187,-35)}

\put(160.4,-20){\circle{2}}
\put(159,-16){$\beta$}
\put(154.6,-26){\circle{2}}
\put(149,-25){$\beta$}
\put(166.8,-24.7){\circle{2}}
\put(166,-21){$\alpha$}
\put(161,-30.8){\circle{2}}
\put(163,-33){$\alpha$}
\put(155,-36.6){\circle{2}}
\put(149,-36){$\alpha$}
\put(160.8,-43){\circle{2}}
\put(159,-48){$\beta$}

\put(60,-60){\path(175,-36)(181,-41)}

\put(60,-60){\path(181,-41)(175,-47)}

\put(60,-60){\path(181,-53)(175,-47)}

\put(235,-95.4){\circle{2}}
\put(237,-93){$\beta$}
\put(240.7,-101){\circle{2}}
\put(243,-101){$\alpha$}
\put(235.4,-106.7){\circle{2}}
\put(230,-106.6){$\alpha$}
\put(240.8,-112.6){\circle{2}}
\put(243,-113){$\beta$}

\put(100,-20){$\path(170,-30)(180,-40)$}

\put(100,-20){\path(180,-40)(175,-45)}

\put(100,-20){\path(180,-50)(175,-45)}

\put(270,-49.3){\circle{2}}
\put(273,-48){$\alpha$}
\put(275.2,-54.3){\circle{2}}
\put(279,-53){$\beta$}
\put(279.5,-59.4){\circle{2}}
\put(283,-59){$\alpha$}
\put(275.3,-65){\circle{2}}
\put(269,-66){$\alpha$}
\put(280,-70.2){\circle{2}}
\put(283,-71){$\beta$}

\put(-20,-45){\path(181,-41)(175,-47)}

\put(-20,-45){\path(181,-53)(175,-47)}

\put(-20,-45){\path(181,-41)(187,-35)}

\put(166.3,-80){\circle{2}}
\put(166,-77){$\alpha$}
\put(160.8,-85.7){\circle{2}}
\put(156,-85){$\alpha$}
\put(155,-91.7){\circle{2}}
\put(149,-91){$\alpha$}
\put(162,-99){\circle{2}}
\put(160,-104){$\beta$}

\put(10,-80){\path(181,-41)(175,-47)}

\put(10,-80){\path(181,-53)(175,-47)}
\put(191,-120.5){\circle{2}}
\put(185,-119){$\alpha$}
\put(184.8,-126.6){\circle{2}}
\put(178,-126){$\alpha$}
\put(190.5,-132.7){\circle{2}}
\put(186.7,-138){$\beta$}

\put(13,-100){\path(181,-53)(175,-47)}

\put(188,-146.8){\circle{2}}
\put(183,-146){$\alpha$}
\put(193.5,-153){\circle{2}}
\put(189,-157.7){$\beta$}

\put(194,-71){\circle{2}}
\put(191,-76.2){$\beta$}
\put(206.7,-70){\circle{2}}
\put(200,-70){$\alpha$}
\put(201,-76){\circle{2}}
\put(201,-80){$\alpha$}
\put(194.7,-82){\circle{2}}
\put(189,-82){$\alpha$}
\put(201,-87.8){\circle{2}}
\put(201.7,-91){$\beta$}

\put(188,-173){$(\beta)$}

\put(188,-200){$(\varnothing)$}

\end{picture}

\n

We will generalize this construction for all irreducible representations of any rank two semisimple Lie algebra.
We want to define the poset $P_{a,b}$ associated to the representation $\Gamma_{a,b}$
in such a way that $L_{a,b}$ gives us the possible paths in $\Gamma_{a,b}$.\\
We need some definitions ( see \cite{[ADLMPPrW]} ).
\begin{defns}

\n
\begin{itemize}
\item[1)] Let $(P, \leq) $ be a partially ordered set and $v,w \in P$ such that $v \leq w$.
We define the interval $[v,w]$ as the set
$$ [v,w]=\{x \in P: v \leq x \leq w\}.$$
We say that $w$ covers $v$ if $[v,w]=\{v,w\}.$\\
\item[2)] A two-color poset is a poset $P$ for which we can associate for each vertex in $P$ a color $\alpha$ or $\beta$.
The function $v \longmapsto color(v)$ is the {\bf{color}} function .\\
\item[3)] We are going to select and numbered some chains in $P$. To do this, we define a chain function:
$$
\hbox{{\bf{chain}}}: P \longrightarrow [[1,m]]
$$
such that:

\begin{itemize}
\item [i)] for $1 \leq i \leq m$ , $chain^{-1}(i)$ \hbox{is~ a ~(possibly~empty)~chain ~in~} $P$.
\item [ii)] $\forall~ u,v \in P$ , if $v$ covers $u$ then either $\hbox{chain}(u)=\hbox{chain}(v)$ or
$\hbox{chain}(u)= \hbox{chain}(v)+1$.
\end{itemize}

\end{itemize}
\end{defns}

\n

\n We represent the function $\hbox{chain}$ as follow:

\n

If $\hbox{chain}(u)= \hbox{chain}(v)+1=k+1 $ then we draw: \vskip-0.5cm
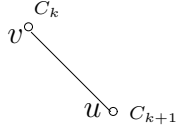
$$\begin{picture}(40,70)(100,-10) \toutpetit

\path(75,32)(105,2)

\put(112,-1){$C_{k+1}$}

\put(74,34){\circle{3}}

\put(106,2){\circle{3}}

\put(76,39){$C_k$}
\end{picture}
\begin{picture}(40,70)(100,-10)
\put(26,28){$v$}

\put(55,0){$u$}
\end{picture}
$$

\n

and if $\hbox{chain}(u)= \hbox{chain}(v)=k$ then we draw: \vskip-0.5cm
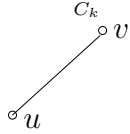
$$\begin{picture}(40,70)(60,-10) \toutpetit

\path(75,32)(42,2)

\put(76,34){\circle{3}}

\put(42,2){\circle{3}}

\put(65,40){$C_k$}
\end{picture}
\begin{picture}(40,70)(100,-10)
\put(80,31){$v$}

\put(46,-3){$u$}
\end{picture}
$$

\begin{exples}

\n

For the $C_2$ case, we shall choose:
\vskip-6cm
$P_{0,1}$:
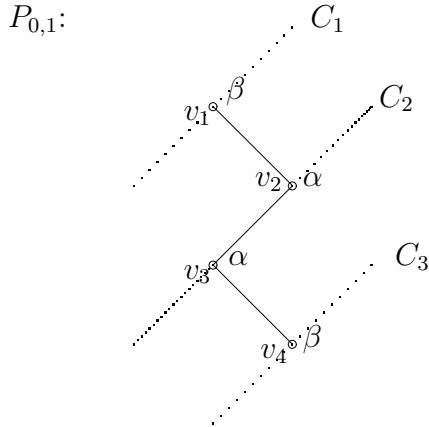
\begin{picture}(100,200)(120,0)\toutpetit

\path(170,-30)(200,-60)

\path(200,-60)(170,-90)

\path(200,-120)(170,-90)

\bezier{20}(200,-60)(230,-30)(230,-30)

\bezier{20}(170,-90)(155,-105)(140,-120)

\bezier{20}(140,-60)(170,-30)(200,0)

\bezier{20}(230,-90)(200,-120)(170,-150)

\put(170,-30){\circle{3}}

\put(200,-60){\circle{3}}

\put(170,-90){\circle{3}}

\put(200,-120){\circle{3}}

\end{picture}
\begin{picture}(40,70)(100,-10)
\put(62,-70){$v_2$}

\put(35,-45){$v_1$}

\put(35,-105){$v_3$}

\put(64,-135){$v_4$}

\put(51,-37){$\beta$}

\put(80,-70){$\alpha$}

\put(52,-100){$\alpha$}

\put(80,-131){$\beta$}

\put(83,-10){$C_1$}

\put(109,-40){$C_2$}

\put(115,-100){$C_3$}

\end{picture}
\vskip6cm
For the $G_2$ case, we choose:
\vskip-2cm

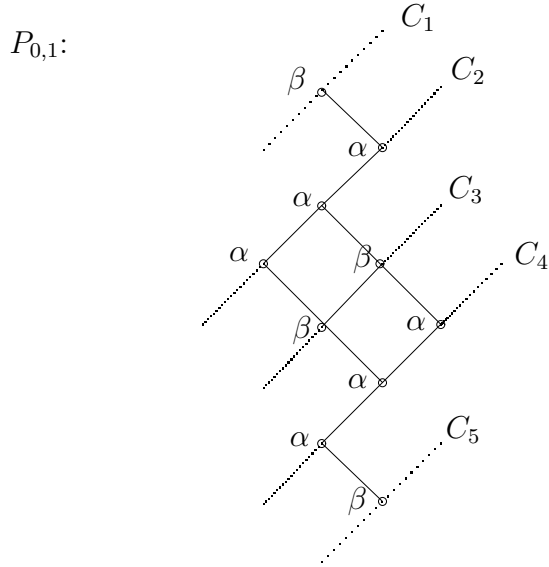
\begin{picture}(100,250)(50,-180)\toutpetit

\path(170,-30)(215,-75)

\path(192,-98)(215,-75)

\path(192,-98)(215,-120)

\path(237,-53)(215,-75)

\path(237,-53)(192,-8)

\path(192,-8)(170,-30)

\path(192,-8)(215,14)

\path(192,36)(215,14)

\path(192,-54)(215,-30)

\bezier{20}(215,58)(192,36)(170,13)

\bezier{20}(237,37)(225,24)(215,14)

\bezier{20}(147,-53)(157,-43)(170,-30)

\bezier{20}(237,-8)(225,-20)(215,-30)

\bezier{20}(192,-54)(180,-67)(170,-77)

\bezier{20}(260,-30)(240,-50)(237,-53)

\bezier{20}(192,-98)(180,-111)(170,-121)

\bezier{20}(237,-98)(215,-120)(192,-143)

\put(170,-30){\circle{3}}

\put(215,-75){\circle{3}}

\put(192,-98){\circle{3}}

\put(215,-120){\circle{3}}

\put(237,-53){\circle{3}}

\put(192,-8){\circle{3}}

\put(192,35){\circle{3}}

\put(215,14){\circle{3}}

\put(192,-54){\circle{3}}

\put(214,-30){\circle{3}}

\end{picture}
\begin{picture}(100,250)(50,-180)

\put(-30,50){$P_{0,1}$:}

\put(75,36){$\beta$}

\put(98,11){$\alpha$}

\put(77,-8){$\alpha$}

\put(53,-29){$\alpha$}

\put(100,-31){$\beta$}

\put(77,-58){$\beta$}

\put(120,-55){$\alpha$}

\put(98,-77){$\alpha$}

\put(76,-99){$\alpha$}

\put(98,-123){$\beta$}

\put(135,-96){$C_5$}

\put(161,-30){$C_4$}

\put(136,-5){$C_3$}

\put(137,40){$C_2$}

\put(118,60){$C_1$}
\end{picture}

\end{exples}

\n

These pictures represent the fundamental posets with the function color and the function chain.
They are uniquely defined with the grid property.
\begin{defn}

\n

A two-color grid poset is a poset $(P, \leq)$ together with a chain function chain and a
color function color such that : if $u$ and $v$ are two vertices in the same connected components
of $P$ and satisfying:
\begin{itemize}

\item [i)]  $ \hbox{if} ~\hbox{chain}(u)= \hbox{chain}(v)+1$ then $ \hbox{color}(u) \neq \hbox{color}(v),$
\item [ii)] $\hbox{if} ~\hbox{chain}(u)= \hbox{chain}(v)$ then $ \hbox{color}(u) = \hbox{color}(v).$
\end{itemize}
\end{defn}
\begin{rem}

\n

On the fundamental posets, there is an unique chain map such that the result are the two-color grid posets.
This choice corresponds to our drawing for each $P_{a,b}$ where $a+b =1$.
\end{rem}
Let us consider now the definition for posets $P_{a,b}$, $a+b \geq 1$.
\begin{defn}

\n

A {\bf{grid}} is a two-color grid poset which has moreover the following max property :
\begin{itemize}
\item [i)] if $u$ is any maximal element in the poset $P$, then
$$
\hbox{chain}(u) \leq \inf_{x \in P} ~~\hbox{chain}(x)+1,
$$
\item [ii)] if $v \neq u$ is another maximal element in $P$, then
$$
color(u) \neq color(v).
$$
\end{itemize}
\end{defn}
\begin{rem}

\n

The fundamental posets are grid posets.
\end{rem}

\n

\n From now one, we identify two grid posets with the same poset, the same color function and two chain maps:
$chain(u)$ and $chain^{'}(u)$, if there is $k$ such that
$chain^{'}(u)= chain(u)+k $ for any $u$.

\begin{defn}

\n
Given two grid posets $P$ and $Q$, we denote by $P \triangleleft Q$ the grid poset with the following properties:
\begin{itemize}

\item [i)]  The elements of $P \triangleleft Q$ is the union of elements of $P$ and those of $Q$.

\item[ii)]  $P$ is an ideal of $P \triangleleft Q$ i.e if $u \in P$ and $v \leq u$ in $P \triangleleft Q$
then $v \in P$, the functions color and chain of $P$ are the restriction of
the functions color and chain of $P \triangleleft Q$ (up to a renumbering of chains ).

\item [iii)] $ ( P \triangleleft Q ) \backslash P $ with the restriction of functions color and chain
on $P \triangleleft Q$ is isomorphic to $Q$ (up to a renumbering of chains p).

\item [iv)] The following holds :
$$
\hbox{if $u$ (resp. $v$) is a maximal element in $P$ (resp. in $Q$), then $chain(u)\leq chain(v)$},
$$
$$
\hbox{ and}
$$
$$
\hbox{if $u$ (resp. $v$) is a minimal element in $P$ (resp. in $Q$), then $chain(u)\leq chain(v)$}.
$$
\end{itemize}
\end{defn}

\n

If  $P$ $\triangleleft$ $Q $ exists, thus $P$ $\triangleleft$ $Q $ is uniquely determined by these conditions,
up to a renumbering of chain.
\begin{rem}

\n

Given three grid posets $P$, $Q$, and $R$ then :
$$
(P \triangleleft Q) \triangleleft R \simeq P \triangleleft (Q \triangleleft R).
$$
We denote this $P \triangleleft Q \triangleleft R$ .
\end{rem}

Starting with the grid posets $P_{1,0}$ and $P_{0,1}$ defined for the rank two semisimple Lie algebra,
for any natural numbers $a$ and $b$, there exists one and only one grid poset
$$ \begin{array}{cccc}
P_{a,b}= & \underbrace{P_{0,1} \triangleleft... \triangleleft P_{0,1}} \triangleleft & \underbrace{ P_{1,0} \triangleleft ... \triangleleft
P_{1,0}}.&\\
& b & a&\\
\end{array}
$$

Now, given the grid poset $P_{a,b}$, we otain a basis of $\Gamma_{a,b}$ by building the corresponding
distributive lattice $L_{a,b}$ of ideals in $P_{a,b}$ and labelling the vertices of $L_{a,b}$ as follows:\\
we start with the heighest weight Young tableaux of shape $\lambda$ : $b$ columns
$\renewcommand{\arraystretch}{0.7}{\begin{array}{l}\framebox{$1$}\\
\framebox{$2$}\end{array}}$ and $a$ columns $\framebox{$1$}$. We put
this tableau on the vertex of $L_{a,b}$ corresponding to the total ideal
$P_{a,b}$. Now, we reach any vertex in $L_{a,b}$ by following a
sequence of edges $\alpha$ or $\beta$. By construction, we know if
this edge corresponds to a vertex in $P_{0,1}$ or in $P_{1,0}$. if
the corresponding vertex is in a $P_{1,0}$-component in $P_{a,b}$,
we act with the edge on the first possible column with size $1$. And
if it is in a $P_{0,1}$-component in $P_{a,b}$, we act with the edge
on the first possible column with size $2$.

\n

Now, we just draw the $L_{2,0}, L_{1,1}$ and  $L_{0,2}$ pictures for each rank two Lie algebra and we call semi
standard tableaux the obtained basis. We summarize the result here:

\begin{prop}

\n

Let  $a$, $b$ be 2 natural numbers, and let $\lambda= (a,b)$. The set of semi standard tableaux for the Lie
algebras of type `type' with size $\lambda$ is denoted $\mathcal{S}_{type}(\lambda)$. Then we get:
$$
\bullet~~\mathcal{S}_{A_1 \times A_1}(\lambda) =\Big\{~ \hbox{usual semi standard tableaux $T$ of shape
$\lambda$ with entries in} ~\{1,2,3,4\}
$$
$$
\hbox{such that}~\framebox{2}~, ~\framebox{4}~, ~\renewcommand{\arraystretch}{0.7}\begin{array}{l}\framebox{1}\\
\framebox{3}\end{array}, \renewcommand{\arraystretch}{0.7}\begin{array}{l}\framebox{2}\\
\framebox{3}\end{array}, \renewcommand{\arraystretch}{0.7}\begin{array}{l}\framebox{2}\\
\framebox{4}\end{array} ~ \hbox{and} ~ \renewcommand{\arraystretch}{0.7}\begin{array}{l}\framebox{3}\\
\framebox{4}\end{array} ~\hbox{are not a column of $T$}\Big\}.
$$

$$
\bullet~~ \mathcal{S}_{A_2}(\lambda) =\Big\{~ \hbox{usual semi standard tableaux $T$ of shape $\lambda$ with
entries in} ~\{1,2,3\}\Big\}.
$$

$$
\bullet~~ \mathcal{S}_{C_2}(\lambda)=\Big\{~ \hbox{usual semi standard tableaux $T$ of shape $\lambda$ with
entries in} ~\{1,2,3,4\}
$$
$$\hbox{such that} ~\renewcommand{\arraystretch}{0.7}\begin{array}{l}\framebox{1}\\
\framebox{4}\end{array} \hbox{is not a column of $T$ and } ~\renewcommand{\arraystretch}{0.7}
\begin{array}{l}\framebox{2}\\
\framebox{3}\end{array} \hbox{appears at most once in $T$}\Big\}.
$$

$$
\bullet~~ \mathcal{S}_{G_2}(\lambda) = \Big\{~ \hbox{usual semi standard tableaux $T$ of shape $\lambda$ with
entries in} ~\{1,2,3,4,5,6,7\}
$$
$$
\hbox{such that the column}~\framebox{4}~\hbox{ appears at most once in $T$,}
\renewcommand{\arraystretch}{0.7} \begin{array}{l}\framebox{2}\\
\framebox{3}\end{array}, \renewcommand{\arraystretch}{0.7}\begin{array}{l}\framebox{2}\\
\framebox{4}\end{array}, \renewcommand{\arraystretch}{0.7}\begin{array}{l}\framebox{3}\\
\framebox{4}\end{array}, \renewcommand{\arraystretch}{0.7}\begin{array}{l}\framebox{3}\\
\framebox{5}\end{array}, \renewcommand{\arraystretch}{0.7}\begin{array}{l}\framebox{4}\\
\framebox{5}\end{array}, \renewcommand{\arraystretch}{0.7}\begin{array}{l}\framebox{4}\\
\framebox{6}\end{array}
$$
$$
\hbox{and} ~ \renewcommand{\arraystretch}{0.7}\begin{array}{l}\framebox{5}\\
\framebox{6}\end{array} ~\hbox{are not a column in $T$ plus the restriction given by the following table}\Big\}.
$$

\end{prop}
\newpage

\vskip1cm

\renewcommand{\arraystretch}{2.5}
\begin{center}
\begin{tabular}{|c|c|}\hline

 Column $T ^{(i)}$ of T

&  Then
the succeeding column $T ^{(i+1)}$ of T cannot be...  \\
\hline
\toutpetit {\framebox{$4$}}  &  \toutpetit {\framebox{$4$}}
\\
\hline
\renewcommand{\arraystretch}{1.4}\toutpetit{$\begin{array}{l}\framebox{$1$}\\
\framebox{$4$}\end{array}$}  &  \renewcommand{\arraystretch}{1.4} \toutpetit {$\framebox{$1$}~,
\begin{array}{l}\framebox{$1$}\\
\framebox{$4$}\end{array} , \begin{array}{l}\framebox{$1$}\\
\framebox{$5$}\end{array}, \begin{array}{l}\framebox{$1$}\\
\framebox{$6$}\end{array}, \begin{array}{l}\framebox{$1$}\\
\framebox{$7$}\end{array}$} 
\\ \hline

\renewcommand{\arraystretch}{1.4}\toutpetit{$\begin{array}{l}\framebox{$1$}\\
\framebox{$5$}\end{array}$} &  \renewcommand{\arraystretch}{1.4}
\toutpetit { $\framebox{1}~,\begin{array}{l}\framebox{$1$}\\
\framebox{5}\end{array}, \begin{array}{l}\framebox{$1$}\\
\framebox{$6$}\end{array}, \begin{array}{l}\framebox{$1$}\\
\framebox{$7$}\end{array}$}
\\ \hline

\renewcommand{\arraystretch}{1.4} \toutpetit  {$\begin{array}{l}\framebox{$1$}\\
\framebox{$6$}\end{array}$}  & \renewcommand{\arraystretch}{1.4}
\toutpetit {$\framebox{$1$}~, \framebox{$2$}~,
 \begin{array}{l}\framebox{$1$}\\
\framebox{$6$}\end{array} , \begin{array}{l}\framebox{$1$}\\
\framebox{7}\end{array}, \begin{array}{l}\framebox{$2$}\\
\framebox{$6$}\end{array}, \begin{array}{l}\framebox{$2$}\\
\framebox{$7$}\end{array}$}
\\ \hline

\renewcommand{\arraystretch}{1.4} \toutpetit  {$\begin{array}{l}\framebox{$2$}\\
\framebox{$6$}\end{array}$}  &  \renewcommand{\arraystretch}{1.4}
\toutpetit {$ \framebox{$2$}~,
\begin{array}{l}\framebox{$2$}\\
\framebox{$6$}\end{array}, \begin{array}{l}\framebox{$2$}\\
\framebox{$7$}\end{array}$} 
\\ \hline

\renewcommand{\arraystretch}{1.4}\toutpetit  {$\begin{array}{l}\framebox{$1$}\\
\framebox{$7$}\end{array}$}  & \renewcommand{\arraystretch}{1.4}
\toutpetit {$\framebox{$1$}~, \framebox{$2$}~, \framebox{$3$}~,
~\framebox{$4$}~,
 \begin{array}{l}\framebox{$1$}\\
\framebox{$7$}\end{array} , \begin{array}{l}\framebox{$2$}\\
\framebox{7}\end{array}, \begin{array}{l}\framebox{$3$}\\
\framebox{$7$}\end{array}, \begin{array}{l}\framebox{$4$}\\
\framebox{$7$}\end{array}$} 
\\ \hline

\renewcommand{\arraystretch}{1.4} \toutpetit  {$\begin{array}{l}\framebox{$2$}\\
\framebox{$7$}\end{array}$}  &  \renewcommand{\arraystretch}{1.4}
\toutpetit {$\framebox{$2$}~, \framebox{$3$}~, \framebox{$4$}~,
\begin{array}{l}\framebox{$2$}\\
\framebox{7}\end{array}, \begin{array}{l}\framebox{$3$}\\
\framebox{$7$}\end{array}, \begin{array}{l}\framebox{$4$}\\
\framebox{$7$}\end{array}$} 
\\ \hline

\renewcommand{\arraystretch}{1.4} \toutpetit {$\begin{array}{l}\framebox{$3$}\\
\framebox{$7$}\end{array}$}  &  \renewcommand{\arraystretch}{1.4}
\toutpetit {$\framebox{$3$}~, \framebox{$4$}~,
 \begin{array}{l}\framebox{$3$}\\
\framebox{$7$}\end{array}, \begin{array}{l}\framebox{$4$}\\
\framebox{$7$}\end{array}$} 
\\ \hline

\renewcommand{\arraystretch}{1.4} \toutpetit {$\begin{array}{l}\framebox{$4$}\\
\framebox{$7$}\end{array}$} & 
\renewcommand{\arraystretch}{1.4}\toutpetit{$
{\framebox{$4$}~, \begin{array}{l}\framebox{$4$}\\
\framebox{$7$}\end{array}}$}
\\ \hline

\end{tabular}
\end{center}

\vskip1.5cm

\section{Shape and reduced shape algebras}

\vskip0.7cm

In section 2, We realize $\mathfrak{g}$ as a space of $n \times n$ matrices ($n=4,3,4,7$),
we consider the subgroup $HN^-$ of matrices in $G$ which are lower triangular, we look for polynomial functions in
the entries of these matrices, the space of these functions is $\mathbb{C}[HN^-]$. We get:\\

\n

$ \underline{\underline{For~A_1 \times A_1}}:$\\

$g \in HN^-$ means
$g=  \left(\begin{smallmatrix}h_1&0&0&0\\0&h_2&0&0\\x&0&h_1^{-1}&0\\ 0&y&0&h_2^{-1}\\
\end{smallmatrix}\right)$ and $\mathbb{C}[HN^-]= \mathbb{C}[h_1,h_2, x,y]$,\\

 \n

 $ \underline{\underline{ For ~A_2}}:$\\

$g \in HN^-$ means
$g=  \left(\begin{smallmatrix}h_1&0&0\\x&h_1^{-1}h_2&0\\ z&y&h_2^{-1}\\
\end{smallmatrix}\right)$ and $\mathbb{C}[HN^-]= \mathbb{C}[h_1,h_2, x,y,z]$,\\

\n

 $ \underline{\underline{For~C_2}}:$\\

$g \in HN^-$ means
$g=  \left(\begin{smallmatrix}h_1&0&0&0\\x& h_2&0&0\\z&u&h_2^{-1}&0\\ y&z&-x&h_1^{-1}\\
\end{smallmatrix}\right)$ and $\mathbb{C}[HN^-]= \mathbb{C}[h_1,h_2, x,y,z,u]$,\\

\n

 $ \underline{\underline{For~G_2}}:$\\

$g \in HN^-$ means
$g=  \left(\begin{smallmatrix}h_1&0&0&0&0&0&0\\x&h_1^{-1}h_2&0&0&0&0&0\\y&a&h_2^{-1}&0&0&0&0\\z&-\sqrt{2}ax+\sqrt{2}y&-\sqrt{2}x&1&0&0&0\\
b&-ax^2+xy-\frac{\sqrt{2}}{2}z&-x^2&\sqrt{2}x&h_1^{-1}&0&0\\ c&axy+\frac{\sqrt{2}}{2}az-y^2&xy+\frac{\sqrt{2}}{2}z& -\sqrt{2}y&-a&h_1h_2^{-1}&0\\
-yb-xc-\frac{z^2}{2}&\frac{\sqrt{2}}{2}axz-ab-\frac{\sqrt{2}}{2}yz-c&\frac{\sqrt{2}}{2}xz-b&-z&-y+ax&-x&h_2
\end{smallmatrix}\right)$ \\ and $\mathbb{C}[HN^-]= \mathbb{C}[h_1,h_2,a,b,c,x,y,z]$.\\

As in the $SL(n)$ case, we can consider these polynomial functions as polynomial functions on the group $G$
which are invariant under multiplication on the right size by $N^+$.
Thus, we get:
$$ \mathbb{C}[HN^-] \simeq \mathbb{C}[G]^{N^+}$$
and the semi standard Young tableaux give a basis for $\mathbb{C}[HN^-]$.
\begin{defn}

\n

The {\bf{shape algebra}} $\mathbb{S}_G$ of $G$ is by definition the
vector space $\mathbb{S}_{G}=\displaystyle
\bigoplus_{a,b}~\Gamma_{a,b}$
 equipped with the multiplication defined by the transposition of the natural maps:
$$
\Gamma_{a+a',b+b'} \hookrightarrow \Gamma_{a,b} \otimes \Gamma_{a',b'}.
$$
\end{defn}

Then by construction, the set of semi standard tableaux forms a basis of the
shape algabra and we get :
$$
\mathbb{C}[G]^{N^{+}} \simeq \mathbb{S}_G \simeq  Sym^\bullet(\mathbb C^,\wedge^2\mathbb~\Big {/}~\mathcal{PL}
$$
where $\mathcal{PL}$ is the ideal generated by all the Pl\"ucker relations (internal or not).
From now one, we consider the restriction of the functions in $\mathbb{S}_G$ to the subgroup $N^-$. We get
a quotient of $\mathbb{S}_G$  which is, as a vector space, the space $\mathbb{C}[N^-]$.

\n

The quotient has the form
$$
\mathbb{C}[HN^-]~\Big/~ <
\renewcommand{\arraystretch}{0.7}{\tiny \begin{array}{l}\framebox{$1$}\\
\framebox{$2$}\end{array}} -1,~
\renewcommand{\arraystretch}{0.5}\tiny{\framebox{1}}~ -1~ >~ \simeq ~
\mathbb{C}[N^-].
$$

\n Since the ideal is $N^+$ invariant, we get a structure of $N^+$
module on this space $\mathbb{C}[N^-]$. This structure is defined
by:
$$
(n^+ .f)(n_1^-)= f(~^t n^+ n_1^-).
$$

\begin{defn}

\n

We call {\bf{reduced shape algebra}} and denote $
\mathbb{S}_G^{red}$ this quotient module, as a vector space, $
\mathbb{S}_G^{red} \simeq \mathbb{C}[N^-] $.
\end{defn}

Starting with the lowest weight vector in any $\Gamma_{a,b} \subset
\mathbb{C}[HN^-]$, which is the Young tableau
$$\renewcommand{\arraystretch}{1.2}{\begin{array}{l}
\begin{array}{|c|c|c|c|c|c|}
\hline  n-1&\hskip0.35cm \cdots\hskip0.35cm&n-1&\hskip0.35cm n \hskip0.38cm&\hskip0.35cm \cdots\hskip0.35cm &\hskip0.35cm n \hskip0.38cm\\
\hline  \end{array}\\ \begin{array}{|c|c|c|}
\hline \hskip0.35cm n \hskip0.38cm  & \hskip0.35cm \cdots\hskip0.35cm &\hskip0.35cm n \hskip0.38cm\\
\hline
\end{array}\\
\end{array}},
$$
acting with $N^+$ , we generate exactly $\Gamma_{a,b}$ thus the
canonical mapping $$ \pi : \mathbb{S}_G \longrightarrow
\mathbb{S}_G^{red}$$ induces a bijective map from $\Gamma_{a,b}$
onto $\pi (\Gamma_{a,b}) $.\\
Now, since the heighest weight vector
$\renewcommand{\arraystretch}{1.2}{\begin{array}{l}
\begin{array}{|c|c|c|c|c|c|}
\hline  1&\hskip0.35cm \cdots\hskip0.35cm&1& 1 &\hskip0.35cm \cdots\hskip0.35cm & 1 \\
\hline  \end{array}\\ \begin{array}{|c|c|c|}
\hline  2   & \hskip0.35cm \cdots\hskip0.35cm & 2 \\
\hline
\end{array}\\
\end{array}}$
is the constant function $1$ in $\mathbb{S}_G^{red}$, the $N^+$
module $\mathbb{S}_G^{red}$ is indecomposable and $\pi
(\Gamma_{a',b'}) \subset \pi (\Gamma_{a,b})$ if $a' \leq a$ and $b'\leq b$.\\

\n

\n Thus we have:
$$
\mathbb{S}_G^{red}= \bigcup _{a,b}~\pi (\Gamma_{a,b})
$$
and
$$
\pi (\Gamma_{a,b})= \bigcup_{\begin{smallmatrix}{a'\leq a\,}\\ {b'\leq b}\end{smallmatrix}}~~\pi (\Gamma_{a',b'}).
$$

We look for a basis for $\mathbb{S}_G^{red}$ which will be well
adapted to this " decomposition" of $\mathbb{C}[N^-] =
\mathbb{S}_G^{red}$.

\vskip1.5cm

\section{Quasi standard Young tableaux}

\vskip0.7cm

Let us give now the definition of quasi standard Young tableaux,
generalizing the $\mathfrak{sl}(3)$ case construction. The set of quasi standard tableaux for the Lie
algebras of type `type' with size $\lambda$ is denoted $\mathcal{QS}_{type}(\lambda)$.We use a case-by-case
argument.\\

$\underline{\underline{\hbox{First~ case}: A_1 \times A_1}}:$\\

\n

\n There is no `external' Pl\"ucker relation in this case, thus we just
suppress the trivial column $\renewcommand{\arraystretch}{0.7}{\begin{array}{l}\framebox{$1$}\\
\framebox{$2$}\end{array}} $ and $
\renewcommand{\arraystretch}{0.5}\framebox{1}$
in the semi standard Young tableaux for $A_1 \times A_1$. Thus we get :
$$\mathcal{QS}_{A_1 \times A_1}(\lambda)=\{ T \in \mathcal{S}_{A_1 \times
A_1}(\lambda), ~\hbox{T without~ any ~trivial ~column}\},
$$
and the picture

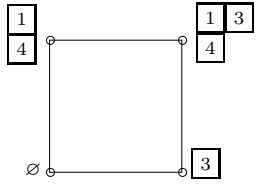
\begin{center}\begin{picture}(150,50)(230,120) \toutpetit

\path(250,150)(300,150)

\path(300,150)(300,100)

\path(250,100)(300,100)

\path(250,150)(250,100)

\put(250.3,150){\circle{3}}

\put(300.3,150){\circle{3}}

\put(300.3,100){\circle{3}}

\put(250.3,100){\circle{3}}

\put(300.5,151){\renewcommand{\arraystretch}{0.7}$\begin{array}{l}
\framebox{1}\framebox{3}\\
\framebox{4}\\
\end{array}$}

\put(229,150.5){\renewcommand{\arraystretch}{0.7}$\begin{array}{l}
\framebox{1}\\
\framebox{4}\\
\end{array}$}

\put(301,101){ \framebox{3}}

\put(238,99){ $\varnothing$}

\end{picture}
\end{center}

\vskip 2cm

$ \underline{\underline{\hbox{Second~ case}:A_2}}:$\\

\n

This case is completely described in section1.\\

$ \underline{\underline{\hbox{Third~ case}:C_2}}:$\\

\n Let $T$ be a semi standard tableau for $C_2$.
 If $T$ does not contain the column
{\renewcommand{\arraystretch}{0.7}$\begin{array}{l}
\framebox{2}\\
\framebox{3}\\
\end{array}$}, we say that $T$ is quasi standard if and only if it
is quasi standard for $\mathfrak{sl}(4)$.  If $T$ contains the column
{\renewcommand{\arraystretch}{0.7}$\begin{array}{l}
\framebox{2}\\
\framebox{3}\\
\end{array}$}, we replace it by the column {\renewcommand{\arraystretch}{0.7}$\begin{array}{l}
\framebox{1}\\
\framebox{4}\\
\end{array}$}, getting a new tableau $T'$.
$T'$ is still semi standard for $\mathfrak{sl}(4)$.\\
\n We say that $T$ is quasi standard if and only if $T'$ is quasi
standard for $\mathfrak{sl}(4)$.
\begin{exple}

\n

For $\lambda= (2,1)$, we get the following family of quasi standard
tableaux with shape $\lambda$:

$ \mathcal{QS}_{C_2}(2,1)=  \left\{\hskip-0.3cm
\begin{array}{llll}\renewcommand{\arraystretch}{0.7}
{\begin{array}{l} \framebox{1}\framebox{3}\framebox{3} \\
\framebox{3} \\
\end{array}},\renewcommand{\arraystretch}{0.7}
{\begin{array}{l}
\framebox{1}\framebox{3}\framebox{4}\\
\framebox{3}\\
\end{array}},\renewcommand{\arraystretch}{0.7}
{\begin{array}{l}
\framebox{1}\framebox{4}\framebox{4}\\
\framebox{3}\\
\end{array}},\renewcommand{\arraystretch}{0.7}
{\begin{array}{l}
\framebox{2}\framebox{4}\framebox{4}\\
\framebox{3}\\
\end{array}},\renewcommand{\arraystretch}{0.7}
{\begin{array}{l}
\framebox{2}\framebox{2}\framebox{2}\\
\framebox{4}\\
\end{array}}
,\renewcommand{\arraystretch}{0.7} {\begin{array}{l}
\framebox{2}\framebox{2}\framebox{3}\\
\framebox{4}\\
\end{array}},\\
\renewcommand{\arraystretch}{0.7} {\begin{array}{l}
\framebox{2}\framebox{3}\framebox{3}\\
\framebox{4}\\
\end{array}},
\renewcommand{\arraystretch}{0.7} {\begin{array}{l}
\framebox{2}\framebox{2}\framebox{4}\\
\framebox{4}\\
\end{array}},\renewcommand{\arraystretch}{0.7}
{\begin{array}{l}
\framebox{2}\framebox{3}\framebox{4}\\
\framebox{4}\\
\end{array}},\renewcommand{\arraystretch}{0.7}
{\begin{array}{l}
\framebox{2}\framebox{4}\framebox{4}\\
\framebox{4}\\
\end{array}},\renewcommand{\arraystretch}{0.7}
{\begin{array}{l}
\framebox{3}\framebox{3}\framebox{3}\\
\framebox{4}\\
\end{array}},\renewcommand{\arraystretch}{0.7}
{\begin{array}{l}
\framebox{3}\framebox{3}\framebox{4}\\
\framebox{4}\\
\end{array}},\\
\renewcommand{\arraystretch}{0.7} {\begin{array}{l}
\framebox{3}\framebox{4}\framebox{4}\\
\framebox{4}\\
\end{array}}
\end{array} \right\}. $

\end{exple}
\begin{thm}

\n

For any $\lambda=(a,b)$, a basis for $\pi(\Gamma_{a,b})$ is
parametrized by the disjoint union
$$
\bigsqcup_{{a'\leq a\,}\atop {b'\leq b}}~~\mathcal{QS}_{C_2} (a',b').
$$

Then, the family of quasi standard Young tableaux forms a basis for
the reduced shape algebra $\mathbb{S}_{C_2}^{red}$.
\end{thm}

{\bf{Proof :}}

\n

We consider $\mathbb{S}_{C_2}^{red}$ as the quotient of the
polynomial algebra in the variables:
$$ X= \framebox{2}~, Y= \framebox{3}~, Z= \framebox{4}~, U= \renewcommand{\arraystretch}{0.7}{\begin{array}{l}
\framebox{1}\\
\framebox{3}\\
\end{array}}~, V= \renewcommand{\arraystretch}{0.7}{\begin{array}{l}
\framebox{2}\\
\framebox{4}\\
\end{array}}~,W=\renewcommand{\arraystretch}{0.7} {\begin{array}{l}
\framebox{2}\\
\framebox{3}\\
\end{array}}~ \hbox{and }~T=\renewcommand{\arraystretch}{0.7} {\begin{array}{l}
\framebox{3}\\
\framebox{4}\\
\end{array}}.$$
By the reduced Pl\"ucker relations : the ideal $\mathcal{PL}_{red}$
generated by the relations:
$$V-XZ+X^2U+Y,~ T-Z^2+ZXU+YU~\hbox{and}~ W+XU-Z.$$

Using the monomial ordering given by the lexicographic ordering on
$(X,Z,Y,W,V,U,T)$, we get the following Groebner basis for
$\mathcal{PL}_{red}$:

Groebner basis of $\mathcal{PL}_{red}$:
$$
\Big\{W^2 +UV-T~,WT+WYU+ZUV-ZT~,-T-YU+ZW,-WY+XT-ZV ,$$
$$W+XU-Z~, -V-Y+XW \Big\}$$
The leading monomials of these elements, with respect to our
ordering are:
$$W^2~, ~ZUV~, ~ZW~, ~XT~, ~XU~,~XW.$$
Thus a basis for the quotient $\mathbb{S}_G^{red}$ is given by the
Young tableaux without any trivial column and which do not contain
the following subtableaux:
$$\renewcommand{\arraystretch}{0.7}\begin{array}{l}\framebox{$2$}\framebox{$2$}\\
\framebox{3}\framebox{$3$}\end{array}~,
\renewcommand{\arraystretch}{0.7}\begin{array}{l}\framebox{$1$}\framebox{$2$}\framebox{$3$}\\
\framebox{3}\framebox{$4$}\end{array}~,
\renewcommand{\arraystretch}{0.7}\begin{array}{l}\framebox{$2$}\framebox{$3$}\\
\framebox{3}\end{array}~,
\renewcommand{\arraystretch}{0.7}\begin{array}{l}\framebox{$3$}\framebox{$2$}\\
\framebox{4}\end{array}~,
\renewcommand{\arraystretch}{0.7}\begin{array}{l}\framebox{$1$}\framebox{$2$}\\
\framebox{3}\end{array}~,
\renewcommand{\arraystretch}{0.7}\begin{array}{l}\framebox{$2$}\framebox{$2$}\\
\framebox{3}\end{array}.
$$
The remaining Young tableaux are exactly the quasi standard Young
tableaux.\\
Indeed, ``$T$ is semi standard without any trivial column'' is
equivalent to ``$T$ does not contain any trival column and does not
contain $\renewcommand{\arraystretch}{0.7}\begin{array}{l}\framebox{$2$}\framebox{$2$}\\
\framebox{3}\framebox{$3$}\end{array}~$ nor
$\renewcommand{\arraystretch}{0.7}
\begin{array}{l}\framebox{$3$}\framebox{$2$}\\
\framebox{4}\end{array}$''. \\
Moreover the remaining tableaux {\sl i.e} $
\renewcommand{\arraystretch}{0.7}\begin{array}{l}\framebox{$1$}\framebox{$2$}\framebox{$3$}\\
\framebox{3}\framebox{$4$}\end{array}~,
\renewcommand{\arraystretch}{0.7}\begin{array}{l}\framebox{$2$}\framebox{$3$}\\
\framebox{3}\end{array}~,\renewcommand{\arraystretch}{0.7}\begin{array}{l}\framebox{$1$}\framebox{$2$}\\
\framebox{3}\end{array}~,\hbox{and}~ \renewcommand{\arraystretch}{0.7}\begin{array}{l}\framebox{$2$}\framebox{$2$}\\
\framebox{3}\end{array}$ are by definition non quasi standards. Now,
if $T$ is a semi standard non quasi standard tableau, without any
trivial column, $T$ contains a minimal semi standard non quasi
standard tableau without trivial column. Looking at all the
possibilities for such minimal tableau with $2$ columns, we get
$$
\renewcommand{\arraystretch}{0.7}\begin{array}{l}\framebox{$1$}\framebox{$2$}\\
\framebox{3}\end{array}~,\renewcommand{\arraystretch}{0.7}\begin{array}{l}\framebox{$2$}\framebox{$2$}\\
\framebox{3}\end{array}~\hbox{and} \renewcommand{\arraystretch}{0.7}\begin{array}{l}\framebox{$2$}\framebox{$3$}\\
\framebox{3}\end{array}.
$$
But there is also such minimal tableau
with three columns. By minimality, such tableau has two columns of
size $2$ and one column of size $1$, $T$ being non quasi standard,
the first column of $T$ is $\renewcommand{\arraystretch}{0.7}
\begin{array}{l}\framebox{$1$}\\
\framebox{3}\end{array}$ or $\renewcommand{\arraystretch}{0.7}
\begin{array}{l}\framebox{$2$}\\
\framebox{3}\end{array}$. \\
If it is
$\renewcommand{\arraystretch}{0.7}
\begin{array}{l}\framebox{$2$}\\
\framebox{$3$}\end{array}$ then we get the non quasi standard
tableaux :
$$
\renewcommand{\arraystretch}{0.7}\begin{array}{l}\framebox{$2$}\framebox{$2$}\framebox{$\makebox(8,8){$u$}$}\\
\framebox{3}\framebox{$4$}\end{array}~, \hbox{and}
\renewcommand{\arraystretch}{0.7}
{\begin{array}{l}\framebox{$2$}\framebox{$3$}\framebox{$\makebox(8,8){$v$}$}\\
\framebox{3}\framebox{$4$}\end{array}}~ \hbox{ with}~ u \geq 2
~\hbox{or}~ v \geq 3.
$$

These non quasi standard tableaux are not
minimal. Thus the first column of $T$ is
$\renewcommand{\arraystretch}{0.7}
\begin{array}{l}\framebox{$1$}\\
\framebox{$3$}\end{array}$ and $T$ is
$$\renewcommand{\arraystretch}{0.7}
\begin{array}{l}\framebox{$1$}\framebox{$2$}\framebox{$\makebox(8,8){$u$}$}\\
\framebox{3}\framebox{$4$}\end{array}~\hbox{ or}
~\renewcommand{\arraystretch}{0.7}
\begin{array}{l}\framebox{$1$}\framebox{$3$}\framebox{$\makebox(8,8){$v$}$}\\
\framebox{3}\framebox{$4$}\end{array}.
$$

The tableau $\renewcommand{\arraystretch}{0.7}
\begin{array}{l}\framebox{$1$}\framebox{$3$}\framebox{$\makebox(8,8){$v$}$}\\
\framebox{3}\framebox{$4$}\end{array}$ are quasi standard for any
$v$. The tableau $\renewcommand{\arraystretch}{0.7}
\begin{array}{l}\framebox{$1$}\framebox{$2$}\framebox{$2$}\\
\framebox{3}\framebox{$4$}\end{array}$ is non quasi standard non
minimal, the tableau $\renewcommand{\arraystretch}{0.7}
\begin{array}{l}\framebox{$1$}\framebox{$2$}\framebox{$3$}\\
\framebox{3}\framebox{$4$}\end{array}$ is non quasi standard
minimal, the tableau $\renewcommand{\arraystretch}{0.7}
\begin{array}{l}\framebox{$1$}\framebox{$2$}\framebox{$3$}\\
\framebox{3}\framebox{$4$}\end{array}$ is quasi standard.\\
The same type of argument shows that any non quasi standard
Young tableaux with more than three columns is not minimal.\\

Finally, if $T$ is any semi standard Young tableau containing a non
quasi standard tableau, $T$ is itself non quasi standard.\\
This proves that the monomial basis for the quotient coincides with
the set of our quasi standard Young tableaux.

\

 $\hfill\square$

Here is the drawing for a part of the diamond cone of $\mathfrak{sp}(4)$\\

{\tiny

\begin{picture}(440,530)(40,0)
 \put(90,180){\vector(0,1){270}}
\put(90,180){\vector(1,-1){180}} \put(137,108){\framebox{2}}
\put(180,48){ \framebox{2}\framebox{2}} \put(211,59){\circle{4}}
\put(150,120){\circle{4}}
\put(83,173){0} \put(91,300){\circle{4}}
\put(91,180){\circle{4}} \put(68,295){$\renewcommand{\arraystretch}{0.7}\begin{array}{l}
\framebox{1}\\
\framebox{3} \\
\end{array}$}

\bezier{60}(150,120)(181,150)(212,180)
\bezier{60}(212,180)(181,210)(150,240)
\bezier{60}(150,240)(120,210)(91,180)
\bezier{20}(91,180)(181,180)(212,180)
\bezier{20}(212,180)(212,240)(212,300)
\bezier{20}(212,300)(181,300)(91,300)

\path(150,120)(273,120)  \put(273,120){\circle{4}}
\put(333,177){\circle{4}} \path(333,177)(333,301)
\path(333,177)(273,120) \path(150,356)(91,301)
\put(150,356){\circle{4}} \path(273,356)(150,356)
\put(273,356){\circle{4}} \path(273,356)(333,301)
\put(334,301){\circle{4}} \put(273,110){$\renewcommand{\arraystretch}{0.7}\begin{array}{l}
\framebox{2}\framebox{2}\\
\framebox{4} \\
\end{array}$}
\put(333,177){$\renewcommand{\arraystretch}{0.7}\begin{array}{l}
\framebox{2}\framebox{4}\\
\framebox{4} \\
\end{array}$}
\put(333,301){$\renewcommand{\arraystretch}{0.7}\begin{array}{l}
\framebox{3}\framebox{4}\\
\framebox{4} \\
\end{array}$}
\put(275,365){$\renewcommand{\arraystretch}{0.7}\begin{array}{l}
\framebox{3}\framebox{3}\\
\framebox{4} \\
\end{array}$}
\put(124,360){$\renewcommand{\arraystretch}{0.7}\begin{array}{l}
\framebox{1}\framebox{3}\\
\framebox{3} \\
\end{array}$}
\put(150,240){\circle{4}} 
\put(212,300){\circle{4}} 
\put(273,240){\circle{4}} 
\put(212,180){\circle{4}} 


\put(201,321){$\renewcommand{\arraystretch}{0.7}\begin{array}{l}
\framebox{2}\framebox{3}\\
\framebox{3} \\
\end{array}$}

\put(215,283){$\renewcommand{\arraystretch}{0.7}\begin{array}{l}
\framebox{3}\\
\framebox{4} \\
\end{array}$}

\put(144,222){\framebox{3}}

\put(140,260){$\renewcommand{\arraystretch}{0.7}\begin{array}{l}
\framebox{2}\\
\framebox{3} \\
\end{array}$}

\put(214,193){$\renewcommand{\arraystretch}{0.7}\begin{array}{l}
\framebox{2}\\
\framebox{4} \\
\end{array}$}

\put(206,160){\framebox{4}}

\put(254,257){$\renewcommand{\arraystretch}{0.7}\begin{array}{l}
\framebox{2}\framebox{4}\\
\framebox{3} \\
\end{array}$}

\put(272,221){$\renewcommand{\arraystretch}{0.7}\begin{array}{l}
\framebox{2}\framebox{3}\\
\framebox{4} \\
\end{array}$}

\end{picture}}

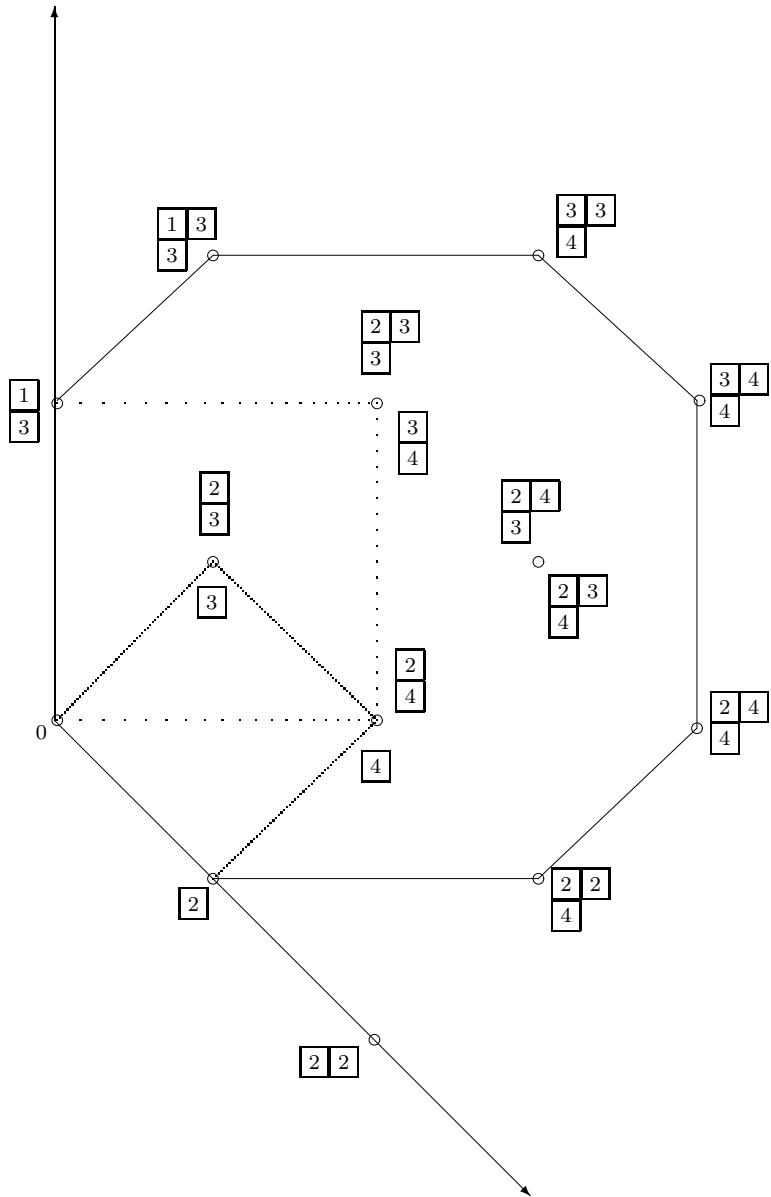
 $\underline{\underline{\hbox{Fourth~ case}:G_2}}$:\\

\begin{defn}

\vskip0.5cm

Let $T= \renewcommand{\arraystretch}{1.2}{\tiny\begin{array}{l}
\begin{array}{|c|c|c|c|c|c|}
\hline \hskip-0.1cm a_1 \hskip-0.15cm &\hskip0.1cm \cdots\hskip0.1cm&\hskip-0.1cm a_p \hskip-0.1cm&\hskip-0.11cm a_{p+1} \hskip-0.11cm &\hskip0.1cm \cdots\hskip0.1cm & \hskip-0.1cm a_{p+q} \hskip-0.1cm \\
\hline  \end{array}\\ \begin{array}{|c|c|c|}
\hline \hskip-0.1cm b_ 1  \hskip-0.12cm & \hskip0.1cm \cdots\hskip0.12cm & \hskip-0.04cm b_p \hskip-0.13cm \\
\hline
\end{array}\\
\end{array}}$ be a semi standard Young tableau of shape $\lambda=(p,q)$ for $G_2$.
We say that $T$ is quasi standard if:
$$
\begin{cases}
       ~{\renewcommand{\arraystretch}{0.7}\tiny{\begin{array}{l}
\framebox{$a_1$\hskip-0.1cm}\\
\framebox{$a_2$\hskip-0.1cm}\\
\end{array}}} \neq {\renewcommand{\arraystretch}{0.7}\tiny{\begin{array}{l}
\framebox{1}\\
\framebox{2}\\
\end{array}}}  \\
~\hbox{and}  \\
       ~ a_1 > 1 \hbox{or ~there~ is~} i=1,...,p ~\hbox{such~ that}~ a_{i+1} > b_i ~\hbox{or}~ a_{i+1} = b_i\neq4. \\
\end{cases}
$$
\end{defn}
Let us denote by $\mathcal{QS}_{G_2}(p,q)$ the set of quasi standard tableaux with shape $(p,q)$, by
$\mathcal{SNQS}_{G_2}(p,q)$ the set of semi standard non quasi standard tableaux with shape $(p,q)$.
We first compute the cardinality of $\mathcal{QS}_{G_2}(p,q)$.

Let us define two operation on $T \in \mathcal{SNQS}_{G_2}(p,q)$.

\n

a) The {\bf{'push' }} operation:\\
Let us denote $T= \renewcommand{\arraystretch}{1.2}{\tiny\begin{array}{l}
\begin{array}{|c|c|c|c|c|c|}
\hline \hskip-0.1cm a_1 \hskip-0.15cm &\hskip0.1cm \cdots\hskip0.1cm&\hskip-0.1cm a_p \hskip-0.1cm&\hskip-0.11cm a_{p+1} \hskip-0.11cm &\hskip0.1cm \cdots\hskip0.1cm & \hskip-0.1cm a_{p+q} \hskip-0.1cm \\
\hline  \end{array}\\ \begin{array}{|c|c|c|}
\hline \hskip-0.1cm b_ 1  \hskip-0.12cm & \hskip0.1cm \cdots\hskip0.12cm & \hskip-0.04cm b_p \hskip-0.13cm \\
\hline
\end{array}\\
\end{array}} \in \mathcal{SNQS}_{G_2}(p,q)$.\\
$\bullet ~~\hbox{If} ~{\renewcommand{\arraystretch}{0.7}\tiny{\begin{array}{l}
\framebox{$a_1$\hskip-0.1cm}\\
\framebox{$a_2$\hskip-0.1cm}\\
\end{array}}}= {\renewcommand{\arraystretch}{0.7}\tiny{\begin{array}{l}
\framebox{1}\\
\framebox{2}\\
\end{array}}}$, we put
$$
P(T)=  \renewcommand{\arraystretch}{1.2}{\tiny\begin{array}{l}
\begin{array}{|c|c|c|c|c|c|}
\hline \hskip-0.1cm a_2 \hskip-0.15cm &\hskip0.1cm \cdots\hskip0.1cm&\hskip-0.1cm a_p \hskip-0.1cm&\hskip-0.11cm a_{p+1} \hskip-0.11cm &\hskip0.1cm \cdots\hskip0.1cm & \hskip-0.1cm a_{p+q} \hskip-0.1cm \\
\hline  \end{array}\\ \begin{array}{|c|c|c|}
\hline \hskip-0.1cm b_ 2  \hskip-0.12cm & \hskip0.1cm \cdots\hskip0.12cm & \hskip-0.04cm b_p \hskip-0.13cm \\
\hline
\end{array}\\
\end{array}}
$$

$ \bullet ~~\hbox{If}~ a_1= 1$ or for any $i=1,...,p$ such that $a_{i+1} < b_i$ or $a_{i+1} = b_i=4$, we put
$$
P(T)=  \renewcommand{\arraystretch}{1.2}{\tiny\begin{array}{l}
\begin{array}{|c|c|c|c|c|c|}
\hline  a_2 &\hskip0.1cm \cdots\hskip0.1cm&\hskip-0.1cm a_{p+1} \hskip-0.15cm&\hskip-0.11cm a_{p+2} \hskip-0.11cm &\hskip0.1cm \cdots\hskip0.1cm & \hskip-0.1cm a_{p+q} \hskip-0.1cm \\
\hline  \end{array}\\ \begin{array}{|c|c|c|}
\hline  b_ 1   & \hskip0.1cm \cdots\hskip0.12cm & b_p \hskip0.1cm \\
\hline
\end{array}\\
\end{array}}.
$$

\n

b) The {\bf{'rectification' }} operation:\\

The tableau $P(T)$ is generally non semi standard. We define the rectification $R(P(T))$ of $P(T)$ as follows:\\
\n we read each $2$ column of $P(T)$ and we replace any wrong $2$ column by a corresponding acceptable one, following
the table 1 :
\vfill
\eject

\begin{center}
\begin{tabular}{|c|c|}\hline

 Wrong column

&  acceptable column  \\
\hline
$\renewcommand{\arraystretch}{1.4}\toutpetit{\begin{array}{l}\framebox{$4$}\\
\framebox{$4$}\end{array}} $ &  $
\renewcommand{\arraystretch}{1.4}\toutpetit {\begin{array}{l}\framebox{$1$}\\
\framebox{$7$}\end{array}}$
\\  \hline
$\renewcommand{\arraystretch}{1.4}\toutpetit{\begin{array}{l}\framebox{$2$}\\
\framebox{$3$}\end{array}}$  & $ \renewcommand{\arraystretch}{1.4} \toutpetit {\begin{array}{l}\framebox{$1$}\\
\framebox{$4$}\end{array}} $
\\ \hline

$\renewcommand{\arraystretch}{1.4}\toutpetit{\begin{array}{l}\framebox{$4$}\\
\framebox{$6$}\end{array}}$ & $ \renewcommand{\arraystretch}{1.4}
\toutpetit {\begin{array}{l}\framebox{$3$}\\
\framebox{7}\end{array}} $
\\ \hline

$ \renewcommand{\arraystretch}{1.4} \toutpetit  {\begin{array}{l}\framebox{$3$}\\
\framebox{$5$}\end{array}}$  & $ \renewcommand{\arraystretch}{1.4}
\toutpetit {\begin{array}{l}\framebox{$2$}\\
\framebox{$6$}\end{array}} $
\\ \hline

$ \renewcommand{\arraystretch}{1.4} \toutpetit  {\begin{array}{l}\framebox{$3$}\\
\framebox{$4$}\end{array}}$  & $ \renewcommand{\arraystretch}{1.4}
\toutpetit {\begin{array}{l}\framebox{$1$}\\
\framebox{$6$}\end{array}} $
\\ \hline

$ \renewcommand{\arraystretch}{1.4}\toutpetit  {\begin{array}{l}\framebox{$5$}\\
\framebox{$6$}\end{array}}$  & $ \renewcommand{\arraystretch}{1.4}
\toutpetit {\begin{array}{l}\framebox{$4$}\\
\framebox{$7$}\end{array}} $
\\ \hline

$ \renewcommand{\arraystretch}{1.4} \toutpetit  {\begin{array}{l}\framebox{$2$}\\
\framebox{$4$}\end{array}}$  & $ \renewcommand{\arraystretch}{1.4}
\toutpetit {\begin{array}{l}\framebox{$1$}\\
\framebox{$5$}\end{array}} $
\\ \hline

$ \renewcommand{\arraystretch}{1.4}\toutpetit  {\begin{array}{l}\framebox{$4$}\\
\framebox{$5$}\end{array}}$  & $ \renewcommand{\arraystretch}{1.4}
\toutpetit {\begin{array}{l}\framebox{$3$}\\
\framebox{$6$}\end{array}} $
\\ \hline

\end{tabular}
\end{center}

\begin{prop}

\n

For any $T \in \mathcal{SNQS}_{G_2}(p,q)$, $R(P(T))$ belongs to
$\mathcal{S}_{G_2}(p-1,q) \sqcup \mathcal{S}_{G_2}(p,q-1)$.

\end{prop}

\n {\bf{Proof:}}

\n

\n If ${\renewcommand{\arraystretch}{0.7}\tiny{\begin{array}{l}
\framebox{$a_1$\hskip-0.1cm}\\
\framebox{$a_2$\hskip-0.1cm}\\
\end{array}}}= {\renewcommand{\arraystretch}{0.7}\tiny{\begin{array}{l}
\framebox{1}\\
\framebox{2}\\
\end{array}}}$ , this is evident. For the second case, using a computer, we consider case by case, all the possibilities
for $3$ successives columns in $T$ and the corresponding result in $P(T)$. We have to consider $3$ cases:\\

$$ \bullet \n \renewcommand{\arraystretch}{1}{\tiny\begin{array}{l}
\begin{array}{|c|c|c|c|cc|}
\hline
\cdots &a_i\hskip-0.12cm &\hskip-0.07cm a_{i+1} \hskip-0.1cm & a_{i+2}\hskip-0.13cm  &  & \cdots   \\
\hline  \end{array} \\ {\tiny\begin{array}{|c|c|c|c|c|}
\hline  \cdots &   b_ i \hskip-0.1cm &  \hskip-0.07cmb_{i+1}\hskip-0.07cm& b_{i+2} \hskip-0.1cm  & \cdots \\
\hline
\end{array}}\\
\end{array}} \stackrel{P}{ \longrightarrow} \renewcommand{\arraystretch}{1}{\tiny\begin{array}{l}
\begin{array}{|c|c|c|cc|}
\hline  \cdots & a_{i+1}\hskip-0.1cm & a_{i+2} \hskip-0.17cm &  &  \cdots    \\
\hline  \end{array} \\ { \tiny\begin{array}{|c|c|c|c|}
\hline  \cdots & \hskip0.1cm   b_ i \hskip0.15cm  &  b_{i+1}\hskip-0.13cm &  \cdots \\
\hline
\end{array}}\\
\end{array}} \stackrel{R}{ \longrightarrow} \renewcommand{\arraystretch}{1}{\tiny\begin{array}{l}
\begin{array}{|c|c|c|cc|}
\hline  \cdots & a'_{i+1}\hskip-0.1cm & a'_{i+2} \hskip-0.17cm &  &  \cdots    \\
\hline  \end{array} \\ { \tiny\begin{array}{|c|c|c|c|}
\hline  \cdots & \hskip0.1cm   b'_ i \hskip0.15cm  &  b'_{i+1}\hskip-0.13cm &  \cdots \\
\hline
\end{array}}\\
\end{array}}$$

$$
\bullet \renewcommand{\arraystretch}{1}{\tiny\begin{array}{l}
\begin{array}{|c|c|c|c|cc|}
\hline
\cdots &a_i\hskip-0.12cm &\hskip-0.07cm a_{i+1} \hskip-0.1cm & a_{i+2}\hskip-0.13cm  &  & \cdots \hskip-0.1cm   \\
\hline  \end{array} \\ {\tiny\begin{array}{|c|c|c|}
\hline  \cdots &   b_ i \hskip-0.1cm &  \hskip-0.07cm b_{i+1}\hskip-0.07cm\\
\hline
\end{array}}\\
\end{array}} \stackrel{P}{ \longrightarrow}\renewcommand{\arraystretch}{1}{\tiny\begin{array}{l}
\begin{array}{|c|c|c|cc|}
\hline  \cdots & a_{i+1}\hskip-0.1cm & a_{i+2} \hskip-0.17cm &  &  \cdots    \\
\hline  \end{array} \\ { \tiny\begin{array}{|c|c|c|}
\hline  \cdots & \hskip0.1cm   b_ i \hskip0.15cm  &  b_{i+1}\hskip-0.13cm \\
\hline
\end{array}}\\
\end{array}} \stackrel{R}{ \longrightarrow} \renewcommand{\arraystretch}{1}{\tiny\begin{array}{l}
\begin{array}{|c|c|c|cc|}
\hline  \cdots & a'_{i+1}\hskip-0.1cm & a'_{i+2} \hskip-0.17cm &  &  \cdots    \\
\hline  \end{array} \\ { \tiny\begin{array}{|c|c|c|}
\hline  \cdots & \hskip0.1cm   b'_ i \hskip0.15cm  &  b'_{i+1}\hskip-0.13cm\\
\hline
\end{array}}\\
\end{array}}
$$

$$
\bullet \renewcommand{\arraystretch}{1}{\tiny\begin{array}{l}
\begin{array}{|c|c|c|c|cc|}
\hline
\cdots &a_i\hskip-0.12cm &\hskip-0.07cm a_{i+1} \hskip-0.1cm & a_{i+2}\hskip-0.13cm  &  & \hskip-0.05cm \cdots \hskip-0.05cm   \\
\hline  \end{array} \\ {\tiny\begin{array}{|c|c|}
\hline  \cdots &  b_ i \hskip-0.1cm \\
\hline
\end{array}}\\
\end{array}} \stackrel{P}{ \longrightarrow} \renewcommand{\arraystretch}{1}{\tiny\begin{array}{l}
\begin{array}{|c|c|c|cc|}
\hline
\cdots &\hskip-0.07cm a_{i+1} \hskip-0.1cm & a_{i+2}\hskip-0.13cm  &  & \hskip-0.05cm \cdots \hskip-0.05cm   \\
\hline  \end{array} \\ {\tiny\begin{array}{|c|c|}
\hline  \cdots &  b_ i \hskip0.15cm \\
\hline
\end{array}}\\
\end{array}} \stackrel{R}{ \longrightarrow} \renewcommand{\arraystretch}{1}{\tiny\begin{array}{l}
\begin{array}{|c|c|c|cc|}
\hline
\cdots &\hskip-0.07cm a'_{i+1} \hskip-0.1cm & a'_{i+2}\hskip-0.13cm  &  & \hskip-0.05cm \cdots \hskip-0.05cm   \\
\hline  \end{array} \\ {\tiny\begin{array}{|c|c|}
\hline  \cdots &  b'_ i \hskip0.15cm \\
\hline
\end{array}}\\
\end{array}}
$$
We verify, in each case, that the result is : $R(P(T)) \subset \mathcal{S}_{G_2}(p,q-1)$.

\n Indeed, for example in the third case, all tableaux $T$ in $\mathcal{S}_{G_2}(1,2)$ such that $ a_{2} < b_1$
define the following tableaux $R(P(T))$:

$$
\left\{\hskip-0.4cm\begin{array}{lll}
\renewcommand{\arraystretch}{0.7}\begin{array}{l}\framebox{$1$}\framebox{$2$}\\
\framebox{3}\end{array},
\renewcommand{\arraystretch}{0.7}\begin{array}{l}\framebox{$1$}\framebox{$1$}\\
\framebox{3}\end{array},
\renewcommand{\arraystretch}{0.7}\begin{array}{l}\framebox{$1$}\framebox{$2$}\\
\framebox{4}\end{array},
\renewcommand{\arraystretch}{0.7}\begin{array}{l}\framebox{$1$}\framebox{$3$}\\
\framebox{3}\end{array},
\renewcommand{\arraystretch}{0.7}\begin{array}{l}\framebox{$1$}\framebox{$3$}\\
\framebox{4}\end{array},
\renewcommand{\arraystretch}{0.7}\begin{array}{l}\framebox{$1$}\framebox{$4$}\\
\framebox{3}\end{array},
\renewcommand{\arraystretch}{0.7}\begin{array}{l}\framebox{$1$}\framebox{$4$}\\
\framebox{4}\end{array},
\renewcommand{\arraystretch}{0.7}\begin{array}{l}\framebox{$1$}\framebox{$5$}\\
\framebox{3}\end{array},
\renewcommand{\arraystretch}{0.7}\begin{array}{l}\framebox{$1$}\framebox{$5$}\\
\framebox{4}\end{array},
\renewcommand{\arraystretch}{0.7}\begin{array}{l}\framebox{$1$}\framebox{$6$}\\
\framebox{3}\end{array},
\renewcommand{\arraystretch}{0.7}\begin{array}{l}\framebox{$1$}\framebox{$6$}\\
\framebox{4}\end{array},\\
\renewcommand{\arraystretch}{0.7}\begin{array}{l}\framebox{$1$}\framebox{$7$}\\
\framebox{3}\end{array},
\renewcommand{\arraystretch}{0.7}\begin{array}{l}\framebox{$1$}\framebox{$7$}\\
\framebox{4}\end{array},
\renewcommand{\arraystretch}{0.7}\begin{array}{l}\framebox{$1$}\framebox{$2$}\\
\framebox{5}\end{array},
\renewcommand{\arraystretch}{0.7}\begin{array}{l}\framebox{$1$}\framebox{$3$}\\
\framebox{5}\end{array},
\renewcommand{\arraystretch}{0.7}\begin{array}{l}\framebox{$1$}\framebox{$3$}\\
\framebox{6}\end{array},
\renewcommand{\arraystretch}{0.7}\begin{array}{l}\framebox{$1$}\framebox{$4$}\\
\framebox{5}\end{array},
\renewcommand{\arraystretch}{0.7}\begin{array}{l}\framebox{$1$}\framebox{$4$}\\
\framebox{6}\end{array},
\renewcommand{\arraystretch}{0.7}\begin{array}{l}\framebox{$1$}\framebox{$5$}\\
\framebox{5}\end{array},
\renewcommand{\arraystretch}{0.7}\begin{array}{l}\framebox{$1$}\framebox{$5$}\\
\framebox{6}\end{array},
\renewcommand{\arraystretch}{0.7}\begin{array}{l}\framebox{$1$}\framebox{$6$}\\
\framebox{5}\end{array},
\renewcommand{\arraystretch}{0.7}\begin{array}{l}\framebox{$1$}\framebox{$6$}\\
\framebox{6}\end{array},\\
\renewcommand{\arraystretch}{0.7}\begin{array}{l}\framebox{$1$}\framebox{$7$}\\
\framebox{5}\end{array},
\renewcommand{\arraystretch}{0.7}\begin{array}{l}\framebox{$1$}\framebox{$7$}\\
\framebox{6}\end{array},
\renewcommand{\arraystretch}{0.7}\begin{array}{l}\framebox{$2$}\framebox{$2$}\\
\framebox{5}\end{array},
\renewcommand{\arraystretch}{0.7}\begin{array}{l}\framebox{$2$}\framebox{$3$}\\
\framebox{5}\end{array},
\renewcommand{\arraystretch}{0.7}\begin{array}{l}\framebox{$2$}\framebox{$3$}\\
\framebox{6}\end{array},
\renewcommand{\arraystretch}{0.7}\begin{array}{l}\framebox{$2$}\framebox{$4$}\\
\framebox{5}\end{array},
\renewcommand{\arraystretch}{0.7}\begin{array}{l}\framebox{$2$}\framebox{$4$}\\
\framebox{6}\end{array},
\renewcommand{\arraystretch}{0.7}\begin{array}{l}\framebox{$2$}\framebox{$5$}\\
\framebox{5}\end{array},
\renewcommand{\arraystretch}{0.7}\begin{array}{l}\framebox{$2$}\framebox{$5$}\\
\framebox{6}\end{array},
\renewcommand{\arraystretch}{0.7}\begin{array}{l}\framebox{$2$}\framebox{$5$}\\
\framebox{7}\end{array},
\renewcommand{\arraystretch}{0.7}\begin{array}{l}\framebox{$2$}\framebox{$6$}\\
\framebox{5}\end{array},\\
\renewcommand{\arraystretch}{0.7}\begin{array}{l}\framebox{$2$}\framebox{$6$}\\
\framebox{6}\end{array},
\renewcommand{\arraystretch}{0.7}\begin{array}{l}\framebox{$2$}\framebox{$6$}\\
\framebox{7}\end{array},
\renewcommand{\arraystretch}{0.7}\begin{array}{l}\framebox{$2$}\framebox{$7$}\\
\framebox{5}\end{array},
\renewcommand{\arraystretch}{0.7}\begin{array}{l}\framebox{$2$}\framebox{$7$}\\
\framebox{6}\end{array},
\renewcommand{\arraystretch}{0.7}\begin{array}{l}\framebox{$2$}\framebox{$7$}\\
\framebox{7}\end{array},
\renewcommand{\arraystretch}{0.7}\begin{array}{l}\framebox{$3$}\framebox{$3$}\\
\framebox{6}\end{array},
\renewcommand{\arraystretch}{0.7}\begin{array}{l}\framebox{$3$}\framebox{$4$}\\
\framebox{6}\end{array},
\renewcommand{\arraystretch}{0.7}\begin{array}{l}\framebox{$3$}\framebox{$5$}\\
\framebox{6}\end{array},
\renewcommand{\arraystretch}{0.7}\begin{array}{l}\framebox{$3$}\framebox{$5$}\\
\framebox{7}\end{array},
\renewcommand{\arraystretch}{0.7}\begin{array}{l}\framebox{$4$}\framebox{$5$}\\
\framebox{7}\end{array},
\renewcommand{\arraystretch}{0.7}\begin{array}{l}\framebox{$3$}\framebox{$6$}\\
\framebox{6}\end{array},\\
\renewcommand{\arraystretch}{0.7}\begin{array}{l}\framebox{$3$}\framebox{$6$}\\
\framebox{7}\end{array},
\renewcommand{\arraystretch}{0.7}\begin{array}{l}\framebox{$4$}\framebox{$6$}\\
\framebox{7}\end{array},
\renewcommand{\arraystretch}{0.7}\begin{array}{l}\framebox{$3$}\framebox{$7$}\\
\framebox{6}\end{array},
\renewcommand{\arraystretch}{0.7}\begin{array}{l}\framebox{$3$}\framebox{$7$}\\
\framebox{7}\end{array},
\renewcommand{\arraystretch}{0.7}\begin{array}{l}\framebox{$4$}\framebox{$7$}\\
\framebox{7}\end{array},
\renewcommand{\arraystretch}{0.7}\begin{array}{l}\framebox{$5$}\framebox{$5$}\\
\framebox{7}\end{array},
\renewcommand{\arraystretch}{0.7}\begin{array}{l}\framebox{$5$}\framebox{$6$}\\
\framebox{7}\end{array},
\renewcommand{\arraystretch}{0.7}\begin{array}{l}\framebox{$6$}\framebox{$6$}\\
\framebox{7}\end{array},
\renewcommand{\arraystretch}{0.7}\begin{array}{l}\framebox{$5$}\framebox{$7$}\\
\framebox{7}\end{array},
\renewcommand{\arraystretch}{0.7}\begin{array}{l}\framebox{$6$}\framebox{$7$}\\
\framebox{7}\end{array},
\renewcommand{\arraystretch}{0.7}\begin{array}{l}\framebox{$1$}\framebox{$2$}\\
\framebox{2}\end{array},\\
\renewcommand{\arraystretch}{0.7}\begin{array}{l}\framebox{$1$}\framebox{$1$}\\
\framebox{2}\end{array}~,
\renewcommand{\arraystretch}{0.7}\begin{array}{l}\framebox{$1$}\framebox{$3$}\\
\framebox{2}\end{array},
\renewcommand{\arraystretch}{0.7}\begin{array}{l}\framebox{$1$}\framebox{$4$}\\
\framebox{2}\end{array},
\renewcommand{\arraystretch}{0.7}\begin{array}{l}\framebox{$1$}\framebox{$5$}\\
\framebox{2}\end{array},
\renewcommand{\arraystretch}{0.7}\begin{array}{l}\framebox{$1$}\framebox{$6$}\\
\framebox{2}\end{array},
\renewcommand{\arraystretch}{0.7}\begin{array}{l}\framebox{$1$}\framebox{$7$}\\
\framebox{2}\end{array}
\end{array} \right\}.
$$
All these tableaux are in $\mathcal{S}_{G_2}(1,1)$.\\

 $\hfill\square$

Now we define a mapping $f$ from $\mathcal{S}_{G_2}(p,q)$ into $\displaystyle\bigsqcup_{{p'\leq p\,}\atop {q' \leq q}}~~\mathcal{QS}_{G_2}(p',q')$
as follows.\\

Let $T$ be in $\mathcal{S}_{G_2}(p,q)$, if $T$ is quasi standard, we put $f(T)=T$, if $T$ is not quasi standard,
we put $T'= R(P(T))$. If $T'$ is quasi standard, we define $f(T)=T'$. If it is not the case, we put $T''= R(P(T'))$,
if $T''$ is quasi standard, we put $f(T)=T''$ and so one...
\begin{prop}

\n

$f$ is a one-to-one onto mapping from $\mathcal{S}_{G_2}(p,q)$ onto
$\displaystyle\bigsqcup_{{p' \leq p\,}\atop {q' \leq q}}~~\mathcal{QS}_{G_2}(p',q')$.
\end{prop}

\n {\bf{Proof:}}

\n

We just define the inverse mapping of $f$. Let $T$ be in $\mathcal{S}_{G_2}(p',q')$. Suppose that $ q' \leq q $.
We first compute $R^{-1}(T)$ i.e we replace each 2-column of $T$ in the ``acceptable columns'' in the table 1 by
the corresponding wrong columns. Let
$$
R^{-1}= \renewcommand{\arraystretch}{1.2}{\tiny\begin{array}{l}
\begin{array}{|c|c|c|c|c|c|}
\hline \hskip-0.1cm a'_1 \hskip-0.15cm &\hskip0.1cm \cdots\hskip0.1cm&\hskip-0.1cm a'_p \hskip-0.1cm&\hskip-0.11cm a'_{p+1} \hskip-0.11cm &\hskip0.1cm \cdots\hskip0.1cm & \hskip-0.1cm a'_{p+q} \hskip-0.1cm \\
\hline  \end{array}\\ \begin{array}{|c|c|c|}
\hline \hskip-0.1cm b'_ 1  \hskip-0.12cm & \hskip0.1cm \cdots\hskip0.12cm & \hskip-0.04cm b'_p \hskip-0.13cm \\
\hline
\end{array}\\
\end{array}}
$$
the resulting tableau. Then we `pull' the resulting tableau, that is we define :
$$
P^{-1}(R^{-1}(T)) = T'=\renewcommand{\arraystretch}{1.2}{\tiny\begin{array}{l}
\begin{array}{|c|c|c|c|c|c|c|}
\hline  1 \hskip0.05cm &a'_1 &\hskip0.1cm \cdots\hskip0.1cm&\hskip-0.1cm a'_{p-1}\hskip-0.07cm \hskip-0.1cm& a'_p  &\hskip0.1cm \cdots\hskip0.1cm & \hskip-0.1cm a'_{p+q} \hskip-0.1cm \\
\hline  \end{array}\\ \begin{array}{|c|c|c|c|}
\hline  b'_ 1 \hskip-0.07cm  & b'_2\hskip0.03cm&\hskip0.1cm \cdots\hskip0.12cm &  b_p \hskip0.07cm\\
\hline
\end{array}\\
\end{array}}.
$$
We verify, case by case as above, that the resulting tableau $T'$ is in $\mathcal{S}_{G_2}(p',q'+1)$.
If $q'+1 < q $, we repeat this operation.\\
Finally, we get a tableau $T''= (P^{-1} \circ R^{-1}) \circ ... \circ (P^{-1} \circ R^{-1})(T) \in
\mathcal{S}_{G_2}(p',q)$. If $ p'< p $, we add to $T''$ $p-p'$ trivial 2-columns
$\renewcommand{\arraystretch}{1.4}\toutpetit{\begin{array}{l}\framebox{$1$}\\
\framebox{$2$}\end{array}}$ . By construction, the mapping $g$ so defined from
$\displaystyle\bigsqcup_{{p' \leq p\,}\atop {q' \leq q}}~~\mathcal{QS}_{G_2}(p',q')$ is the inverse mapping of $f$.\\
Let us recall the projection mapping $ \pi: \mathbb{S}_{G_2}= \oplus_{p,q}~~\Gamma_{q,p} \longrightarrow
\mathbb{S}_{G_2}^{red}$. We show that if $p' \leq p, q' \leq q$, then $\pi(\Gamma_{q',p'}) \subset \Gamma_{q,p}$.
Now, our proposition proves by induction on $p$ and $q$ that:
$$
\sharp \mathcal{QS}_{G_2}(p,q)= dim\big ( \pi (\Gamma_{q,p})~\Big{/}
\displaystyle\sum _{(p',q')<(p,q)}~ \pi (\Gamma_{q',p'})\big)
$$
where $(p',q')<(p,q)$ means $p'\leq p$, $q'\leq q$ and $(p',q')\neq(p,q)$.

\begin{prop}

\n

The set $\mathcal{QS}_{G_2}(p,q)$ is a basis for a supplementary space $W_{p,q}$ in $ \pi(\Gamma_{q,p})$ to
the space $\displaystyle\sum _{(p',q')<(p,q)}~ \pi (\Gamma_{q',p'})$.
\end{prop}

\n {\bf{Proof:}}

\n

Since the number of quasi standard tableaux is the dimension of our space, it is enough to prove that the family
$\mathcal{QS}_{G_2}(p,q)$ is independant in the quotient \\ $\pi (\Gamma_{q,p})~\Big/
\displaystyle\sum _{(p',q')<(p,q)}~ \pi (\Gamma_{q',p'})$.\\

Suppose this is not the case, there is a linear relation $ \displaystyle\sum _{i}~ a_i T_i$  between some
$T_i$ in $\mathcal{QS}_{G_2}(p,q)$ which belongs to $\displaystyle\sum _{(p',q')<(p,q)}~ \pi (\Gamma_{q',p'})$
that means, there is a $S$ in the ideal $\mathcal{PL}_{red}$ of reduced
Pl\"ucker ideal, a family $(T'_j)$ of tableaux in $\cup_{(p',q') < (p,q)}~~\mathcal{S}_{G_2}(p',q')$ and $b_j \in \mathbb{R}$
such that: $\displaystyle\sum _{i}~ a_i T_i= \displaystyle\sum _{j}~ b_j T'_j + S$. This means
$$
\big(\displaystyle\sum _{i}~ a_i T_i- \displaystyle\sum _{j}~ b_j T'_j\big)~|_{N^-}=0.\hskip3cm (1)$$
But now the action of the diagonal matrices $H \in \mathfrak{h}$ in $G_2$ are diagonalized in
$\mathbb{C}[{\mathcal{C}}]$ where ${\mathcal{C}}$ is the set of the polynomial functions
$
\renewcommand{\arraystretch}{0.7}\begin{array}{l}\framebox{$i_1$}\\
\framebox{$i_2$}\\
\framebox{~\vdots\hskip0.055cm} \\
\framebox{$i_r$\hskip0.007cm}\\
\end{array}$ already defined in section 1. Thus we decompose the preceding expression in a finite sum of weight vectors whith weight
$\mu \in \mathfrak{h}^*$. The relation (1) holds for any weight vector, thus we get:
$$
\big(\displaystyle\sum _{i}~ a_i T_i- \displaystyle\sum _{j}~ b_j T'_j\big)~|_{N^-}=0,
$$
$$
H. \big(\displaystyle\sum _{i}~ a_i T_i- \displaystyle\sum _{j}~ b_j T'_j\big) =
\mu(H) \big(\displaystyle\sum _{i}~ a_i T_i- \displaystyle\sum _{j}~ b_j T'_j\big).
$$

\n The first relation means there is $S_{\mu}$ in the ideal $\mathcal{PL}_{red}$ such that:
$$
\displaystyle\sum _{i}~ a_i T_i- \displaystyle\sum _{j}~ b_j T'_j = S_{\mu}.
$$
$S_{\mu}$ being in $\mathcal{PL}_{red}$ can be written as:
$$
S_{\mu}= \displaystyle\sum _{k}~  PL_k + \displaystyle\sum _{l}~ T'_l\big(\renewcommand{\arraystretch}{1.4}
\toutpetit {\begin{array}{l}\framebox{$1$}\\
\framebox{$2$}\end{array}}~-1~\big) + \displaystyle\sum _{m}~ T''_m\big(\renewcommand{\arraystretch}{1.4}
\toutpetit {\begin{array}{l}\framebox{$1$}\end{array}}~-1~\big)
$$
where $PL_k$ are Pl\"ucker relations which are homogeneous, with weight $\mu$, with respect to the $\mathfrak{h}$ action. Let us put
$$
U=\displaystyle\sum _{i}~ a_i T_i- \displaystyle\sum _{j}~ b_j T'_j - \displaystyle\sum _{k}~ PL_k.
$$

$U$ is a linear combination of Young tableaux $U=\sum_\ell c_\ell U_\ell$, it is homogeneous with weight $\mu$.
If we supress the trivial columns of each the $U_\ell$ tableau, we get a tableau $U'_\ell$ of weight
$\mu-a\omega_1-b\omega_2$, if there is $a$ columns $\framebox{$1$}$ and $b$ columns
$\renewcommand{\arraystretch}{0.7}\begin{array}{l}\framebox{$1$}\\
\framebox{2}\end{array}$. Now to supress these columns corresponds exactly to the restriction of the corresponding
polynomial functions to $N^-$. Denoting by $'$ the restriction to $N^-$, we get :
$$
U'=\sum_\ell c_\ell U'_\ell=0.
$$
For any $(a,b)$, we put $M_{(a,b)}=\{\ell,~\hbox{such that $U'_\ell$ has weight }\mu-a\omega_1-b\omega_2\}$ then
for any $(a,b)$, by homogeneity,
$$
\sum_{\ell\in M_{(a,b)}}c_\ell U_\ell=0.
$$

Finally,
$$
U=\sum_{a,b}\left(\renewcommand{\arraystretch}{0.7}\begin{array}{l}\framebox{$1$}\\
\framebox{2}\end{array}\right)^b\sum_{\ell\in M_{(a,b)}}c_\ell U_\ell\left(\framebox{$1$}\right)^a=0.
$$
This proves our proposition.

$\hfill\square$

Finally we can compute the semi standard non quasi standard minimal tableaux for $G_2$, without any trivial
column :
$$ \left\{\hskip-0.3cm\begin{array}{lll}
\renewcommand{\arraystretch}{0.7}\begin{array}{l}\framebox{$1$}\framebox{$2$}\\
\framebox{3}\end{array},
\renewcommand{\arraystretch}{0.7}\begin{array}{l}\framebox{$1$}\framebox{$2$}\\
\framebox{4}\end{array},
\renewcommand{\arraystretch}{0.7}\begin{array}{l}\framebox{$1$}\framebox{$2$}\\
\framebox{5}\end{array},
\renewcommand{\arraystretch}{0.7}\begin{array}{l}\framebox{$1$}\framebox{$3$}\\
\framebox{4}\end{array},
\renewcommand{\arraystretch}{0.7}\begin{array}{l}\framebox{$1$}\framebox{$3$}\\
\framebox{5}\end{array},
\renewcommand{\arraystretch}{0.7}\begin{array}{l}\framebox{$1$}\framebox{$3$}\\
\framebox{6}\end{array},
\renewcommand{\arraystretch}{0.7}\begin{array}{l}\framebox{$1$}\framebox{$4$}\\
\framebox{4}\end{array},
\renewcommand{\arraystretch}{0.7}\begin{array}{l}\framebox{$1$}\framebox{$4$}\\
\framebox{5}\end{array},
\renewcommand{\arraystretch}{0.7}\begin{array}{l}\framebox{$1$}\framebox{$4$}\\
\framebox{6}\end{array},
\renewcommand{\arraystretch}{0.7}\begin{array}{l}\framebox{$1$}\framebox{$5$}\\
\framebox{6}\end{array},\\
\renewcommand{\arraystretch}{0.7}\begin{array}{l}\framebox{$1$}\framebox{$5$}\\
\framebox{7}\end{array},
\renewcommand{\arraystretch}{0.7}\begin{array}{l}\framebox{$1$}\framebox{$6$}\\
\framebox{7}\end{array},
\renewcommand{\arraystretch}{0.7}\begin{array}{l}\framebox{$1$}\framebox{$2$}\framebox{$3$}\\
\framebox{3}\framebox{$5$}\end{array},
\renewcommand{\arraystretch}{0.7}\begin{array}{l}\framebox{$1$}\framebox{$2$}\framebox{$4$}\\
\framebox{3}\framebox{$5$}\end{array},
\renewcommand{\arraystretch}{0.7}\begin{array}{l}\framebox{$1$}\framebox{$2$}\framebox{$3$}\\
\framebox{3}\framebox{$6$}\end{array},
\renewcommand{\arraystretch}{0.7}\begin{array}{l}\framebox{$1$}\framebox{$2$}\framebox{$4$}\\
\framebox{3}\framebox{$6$}\end{array},
\renewcommand{\arraystretch}{0.7}\begin{array}{l}\framebox{$1$}\framebox{$2$}\framebox{$5$}\\
\framebox{3}\framebox{$3$}\end{array},
\renewcommand{\arraystretch}{0.7}\begin{array}{l}\framebox{$1$}\framebox{$2$}\framebox{$5$}\\
\framebox{3}\framebox{$7$}\end{array},\\
\renewcommand{\arraystretch}{0.7}\begin{array}{l}\framebox{$1$}\framebox{$2$}\framebox{$6$}\\
\framebox{3}\framebox{$7$}\end{array},
\renewcommand{\arraystretch}{0.7}\begin{array}{l}\framebox{$1$}\framebox{$2$}\framebox{$4$}\\
\framebox{4}\framebox{$5$}\end{array},
\renewcommand{\arraystretch}{0.7}\begin{array}{l}\framebox{$1$}\framebox{$2$}\framebox{$4$}\\
\framebox{4}\framebox{$6$}\end{array},
\renewcommand{\arraystretch}{0.7}\begin{array}{l}\framebox{$1$}\framebox{$2$}\framebox{$5$}\\
\framebox{4}\framebox{$6$}\end{array},
\renewcommand{\arraystretch}{0.7}\begin{array}{l}\framebox{$1$}\framebox{$3$}\framebox{$4$}\\
\framebox{4}\framebox{$6$}\end{array},
\renewcommand{\arraystretch}{0.7}\begin{array}{l}\framebox{$1$}\framebox{$3$}\framebox{$5$}\\
\framebox{4}\framebox{$6$}\end{array},
\renewcommand{\arraystretch}{0.7}\begin{array}{l}\framebox{$1$}\framebox{$2$}\framebox{$5$}\\
\framebox{4}\framebox{$7$}\end{array},
\renewcommand{\arraystretch}{0.7}\begin{array}{l}\framebox{$1$}\framebox{$2$}\framebox{$6$}\\
\framebox{4}\framebox{$7$}\end{array},\\
\renewcommand{\arraystretch}{0.7}\begin{array}{l}\framebox{$1$}\framebox{$3$}\framebox{$5$}\\
\framebox{4}\framebox{$7$}\end{array},
\renewcommand{\arraystretch}{0.7}\begin{array}{l}\framebox{$1$}\framebox{$3$}\framebox{$6$}\\
\framebox{4}\framebox{$7$}\end{array},
\renewcommand{\arraystretch}{0.7}\begin{array}{l}\framebox{$1$}\framebox{$2$}\framebox{$5$}\\
\framebox{5}\framebox{$6$}\end{array},
\renewcommand{\arraystretch}{0.7}\begin{array}{l}\framebox{$1$}\framebox{$3$}\framebox{$5$}\\
\framebox{5}\framebox{$6$}\end{array},
\renewcommand{\arraystretch}{0.7}\begin{array}{l}\framebox{$1$}\framebox{$2$}\framebox{$5$}\\
\framebox{5}\framebox{$7$}\end{array},
\renewcommand{\arraystretch}{0.7}\begin{array}{l}\framebox{$1$}\framebox{$2$}\framebox{$6$}\\
\framebox{5}\framebox{$7$}\end{array},
\renewcommand{\arraystretch}{0.7}\begin{array}{l}\framebox{$1$}\framebox{$3$}\framebox{$5$}\\
\framebox{5}\framebox{$7$}\end{array},
\renewcommand{\arraystretch}{0.7}\begin{array}{l}\framebox{$1$}\framebox{$3$}\framebox{$6$}\\
\framebox{5}\framebox{$7$}\end{array},\\
\renewcommand{\arraystretch}{0.7}\begin{array}{l}\framebox{$1$}\framebox{$4$}\framebox{$5$}\\
\framebox{5}\framebox{$7$}\end{array},
\renewcommand{\arraystretch}{0.7}\begin{array}{l}\framebox{$1$}\framebox{$4$}\framebox{$6$}\\
\framebox{5}\framebox{$7$}\end{array},
\renewcommand{\arraystretch}{0.7}\begin{array}{l}\framebox{$1$}\framebox{$3$}\framebox{$6$}\\
\framebox{6}\framebox{$7$}\end{array},
\renewcommand{\arraystretch}{0.7}\begin{array}{l}\framebox{$1$}\framebox{$4$}\framebox{$6$}\\
\framebox{6}\framebox{$7$}\end{array},
\renewcommand{\arraystretch}{0.7}\begin{array}{l}\framebox{$1$}\framebox{$5$}\framebox{$6$}\\
\framebox{6}\framebox{$7$}\end{array},
\renewcommand{\arraystretch}{0.7}\begin{array}{l}\framebox{$1$}\framebox{$2$}\framebox{$3$}\framebox{$5$}\\
\framebox{3}\framebox{$5$}\framebox{$6$}\end{array},
\renewcommand{\arraystretch}{0.7}\begin{array}{l}\framebox{$1$}\framebox{$2$}\framebox{$3$}\framebox{$6$}\\
\framebox{3}\framebox{$6$}\framebox{$7$}\end{array},\\
\renewcommand{\arraystretch}{0.7}\begin{array}{l}\framebox{$1$}\framebox{$2$}\framebox{$4$}\framebox{$6$}\\
\framebox{3}\framebox{$6$}\framebox{$7$}\end{array},
\renewcommand{\arraystretch}{0.7}\begin{array}{l}\framebox{$1$}\framebox{$2$}\framebox{$5$}\framebox{$6$}\\
\framebox{3}\framebox{$6$}\framebox{$7$}\end{array},
\renewcommand{\arraystretch}{0.7}\begin{array}{l}\framebox{$1$}\framebox{$2$}\framebox{$3$}\framebox{$5$}\\
\framebox{3}\framebox{$5$}\framebox{$7$}\end{array},
\renewcommand{\arraystretch}{0.7}\begin{array}{l}\framebox{$1$}\framebox{$2$}\framebox{$4$}\framebox{$5$}\\
\framebox{3}\framebox{$5$}\framebox{$7$}\end{array},
\renewcommand{\arraystretch}{0.7}\begin{array}{l}\framebox{$1$}\framebox{$2$}\framebox{$3$}\framebox{$6$}\\
\framebox{3}\framebox{$5$}\framebox{$7$}\end{array},
\renewcommand{\arraystretch}{0.7}\begin{array}{l}\framebox{$1$}\framebox{$2$}\framebox{$4$}\framebox{$6$}\\
\framebox{3}\framebox{$5$}\framebox{$7$}\end{array},\\
\renewcommand{\arraystretch}{0.7}\begin{array}{l}\framebox{$1$}\framebox{$2$}\framebox{$4$}\framebox{$5$}\\
\framebox{4}\framebox{$5$}\framebox{$7$}\end{array},
\renewcommand{\arraystretch}{0.7}\begin{array}{l}\framebox{$1$}\framebox{$2$}\framebox{$4$}\framebox{$6$}\\
\framebox{4}\framebox{$5$}\framebox{$7$}\end{array},
\renewcommand{\arraystretch}{0.7}\begin{array}{l}\framebox{$1$}\framebox{$2$}\framebox{$4$}\framebox{$7$}\\
\framebox{4}\framebox{$5$}\framebox{$7$}\end{array},
\renewcommand{\arraystretch}{0.7}\begin{array}{l}\framebox{$1$}\framebox{$2$}\framebox{$4$}\framebox{$6$}\\
\framebox{4}\framebox{$6$}\framebox{$7$}\end{array},
\renewcommand{\arraystretch}{0.7}\begin{array}{l}\framebox{$1$}\framebox{$2$}\framebox{$5$}\framebox{$6$}\\
\framebox{4}\framebox{$6$}\framebox{$7$}\end{array},
\renewcommand{\arraystretch}{0.7}\begin{array}{l}\framebox{$1$}\framebox{$3$}\framebox{$4$}\framebox{$6$}\\
\framebox{4}\framebox{$6$}\framebox{$7$}\end{array},\\
\renewcommand{\arraystretch}{0.7}\begin{array}{l}\framebox{$1$}\framebox{$3$}\framebox{$5$}\framebox{$6$}\\
\framebox{4}\framebox{$6$}\framebox{$7$}\end{array},
\renewcommand{\arraystretch}{0.7}\begin{array}{l}\framebox{$1$}\framebox{$2$}\framebox{$5$}\framebox{$6$}\\
\framebox{5}\framebox{$6$}\framebox{$7$}\end{array},
\renewcommand{\arraystretch}{0.7}\begin{array}{l}\framebox{$1$}\framebox{$3$}\framebox{$5$}\framebox{$6$}\\
\framebox{5}\framebox{$6$}\framebox{$7$}\end{array},
\renewcommand{\arraystretch}{0.7}\begin{array}{l}\framebox{$1$}\framebox{$2$}\framebox{$3$}\framebox{$5$}\framebox{6}\\
\framebox{3}\framebox{$5$}\framebox{$6$}\framebox{7}\end{array}\end{array}\right\}.$$

\vskip 2cm

Now, for $G_2$, the picture of a part of the diamond cone is as follows :

\newpage

\begin{picture}(100,65)(185,500) \toutpetit

\path(330,-70)(445,-70)

\path(613,35)(445,-70)

\path(613,35)(673,130)

\path(673,130)(673,325)

\path(673,325)(613,420)

\path(613,420)(445,515)

\path(445,515)(330,515)

\path(330,515)(165,420)

\path(165,420)(103,325)

\path(103,325)(103,130)

\path(103,130)(163,35)

\path(330,-70)(163,35)

\put(330,-70){\circle{5}}

\put(445,-70){\circle{5}}

\put(613,35){\circle{5}}

\put(673,130){\circle{5}}

\put(673,325){\circle{5}}

\put(613,420){\circle{5}}

\put(445,515){\circle{5}}

\put(330,515){\circle{5}}

\put(165,420){\circle{5}}

\put(103,325){\circle{5}}

\put(103,130){\circle{5}}

\put(163,35){\circle{5}}

\put(325,-81) {0}

\put(618,33) {$\renewcommand{\arraystretch}{0.7}\begin{array}{l}
\framebox{2}\framebox{2}\\
\framebox{5} \\
\end{array}$}

\put(640,129) {$\renewcommand{\arraystretch}{0.7}\begin{array}{l}
\framebox{2}\framebox{5}\\
\framebox{5} \\
\end{array}$}

\put(640,319) {$\renewcommand{\arraystretch}{0.7}\begin{array}{l}
\framebox{5}\framebox{5}\\
\framebox{7} \\
\end{array}$}

\put(615,424) {$\renewcommand{\arraystretch}{0.7}\begin{array}{l}
\framebox{5}\framebox{7}\\
\framebox{7} \\
\end{array}$}

\put(447,525) {$\renewcommand{\arraystretch}{0.7}\begin{array}{l}
\framebox{6}\framebox{7}\\
\framebox{7} \\
\end{array}$}

\put(308,527) {$\renewcommand{\arraystretch}{0.7}\begin{array}{l}
\framebox{6}\framebox{6}\\
\framebox{7} \\
\end{array}$}

\put(132,422) {$\renewcommand{\arraystretch}{0.7}\begin{array}{l}
\framebox{3}\framebox{6}\\
\framebox{6} \\
\end{array}$}

\put(107,319) {$\renewcommand{\arraystretch}{0.7}\begin{array}{l}
\framebox{3}\framebox{3}\\
\framebox{6} \\
\end{array}$}

\put(107,130) {$\renewcommand{\arraystretch}{0.7}\begin{array}{l}
\framebox{1}\framebox{3}\\
\framebox{3} \\
\end{array}$}

\put(556,130){\circle{5}}

\put(529,138) {$\renewcommand{\arraystretch}{0.7}\begin{array}{l}
\framebox{1}\framebox{5}\\
\framebox{5} \\
\end{array}$}

\put(558,122) {$\renewcommand{\arraystretch}{0.7}\begin{array}{l}
\framebox{2}\framebox{4}\\
\framebox{5} \\
\end{array}$}

\put(612,230){\circle{5}}

\put(613,222) {$\renewcommand{\arraystretch}{0.7}\begin{array}{l}
\framebox{2}\framebox{7}\\
\framebox{5} \\
\end{array}$}

\put(584,237) {$\renewcommand{\arraystretch}{0.7}\begin{array}{l}
\framebox{2}\framebox{5}\\
\framebox{7} \\
\end{array}$}

\put(552,325){\circle{5}}

\put(554,315) {$\renewcommand{\arraystretch}{0.7}\begin{array}{l}
\framebox{2}\framebox{7}\\
\framebox{7} \\
\end{array}$}

\put(520,330) {$\renewcommand{\arraystretch}{0.7}\begin{array}{l}
\framebox{4}\framebox{5}\\
\framebox{7} \\
\end{array}$}

\put(502,420){\circle{5}}

\put(504,410) {$\renewcommand{\arraystretch}{0.7}\begin{array}{l}
\framebox{5}\framebox{6}\\
\framebox{7} \\
\end{array}$}

\put(470,425) {$\renewcommand{\arraystretch}{0.7}\begin{array}{l}
\framebox{4}\framebox{7}\\
\framebox{7} \\
\end{array}$}

\put(387,420){\circle{5}}

\put(389,410) {$\renewcommand{\arraystretch}{0.7}\begin{array}{l}
\framebox{4}\framebox{6}\\
\framebox{7} \\
\end{array}$}

\put(352,420) {$\renewcommand{\arraystretch}{0.7}\begin{array}{l}
\framebox{3}\framebox{7}\\
\framebox{7} \\
\end{array}$}

\put(273,420){\circle{5}}

\put(275,410) {$\renewcommand{\arraystretch}{0.7}\begin{array}{l}
\framebox{3}\framebox{6}\\
\framebox{7} \\
\end{array}$}

\put(238,420) {$\renewcommand{\arraystretch}{0.7}\begin{array}{l}
\framebox{3}\framebox{7}\\
\framebox{6} \\
\end{array}$}

\put(220,325){\circle{5}}

\put(222,315) {$\renewcommand{\arraystretch}{0.7}\begin{array}{l}
\framebox{3}\framebox{4}\\
\framebox{6} \\
\end{array}$}

\put(185,325) {$\renewcommand{\arraystretch}{0.7}\begin{array}{l}
\framebox{1}\framebox{6}\\
\framebox{6} \\
\end{array}$}

\put(442,325){\circle{5}}

\put(440,305) {$\renewcommand{\arraystretch}{0.7}\begin{array}{l}
\framebox{3}\framebox{5}\\
\framebox{7} \\
\end{array}$}

\put(445,331) {$\renewcommand{\arraystretch}{0.7}\begin{array}{l}
\framebox{2}\framebox{7}\\
\framebox{6} \\
\end{array}$}

\put(413,308) {$\renewcommand{\arraystretch}{0.7}\begin{array}{l}
\framebox{2}\framebox{6}\\
\framebox{7} \\
\end{array}$}

\put(420,333) {$\renewcommand{\arraystretch}{0.7}\begin{array}{l}
\framebox{1}\framebox{7}\\
\framebox{7} \\
\end{array}$}

\put(398,230){\circle{5}}

\put(400,220) {$\renewcommand{\arraystretch}{0.7}\begin{array}{l}
\framebox{1}\framebox{7}\\
\framebox{4} \\
\end{array}$}

\put(373,206) {$\renewcommand{\arraystretch}{0.7}\begin{array}{l}
\framebox{1}\framebox{6}\\
\framebox{5} \\
\end{array}$}

\put(365,232) {$\renewcommand{\arraystretch}{0.7}\begin{array}{l}
\framebox{2}\framebox{4}\\
\framebox{6} \\
\end{array}$}

\put(395,248) {$\renewcommand{\arraystretch}{0.7}\begin{array}{l}
\framebox{5}\\
\framebox{6} \\
\end{array}$}

\put(273,230){\circle{5}}

\put(274,218) {$\renewcommand{\arraystretch}{0.7}\begin{array}{l}
\framebox{1}\framebox{6}\\
\framebox{4} \\
\end{array}$}

\put(245,235) {$\renewcommand{\arraystretch}{0.7}\begin{array}{l}
\framebox{1}\framebox{7}\\
\framebox{3} \\
\end{array}$}

\put(248,207) {$\renewcommand{\arraystretch}{0.7}\begin{array}{l}
\framebox{2}\framebox{3}\\
\framebox{6} \\
\end{array}$}

\put(271,245) {$\renewcommand{\arraystretch}{0.7}\begin{array}{l}
\framebox{3}\\
\framebox{7} \\
\end{array}$}

\put(220,130){\circle{5}}

\put(221,125) {$\renewcommand{\arraystretch}{0.7}\begin{array}{l}
\framebox{1}\framebox{3}\\
\framebox{4} \\
\end{array}$}

\put(198,135) {$\renewcommand{\arraystretch}{0.7}\begin{array}{l}
\framebox{1}\\
\framebox{6} \\
\end{array}$}

\put(388,35){\circle{5}}

\put(390,35) {$\renewcommand{\arraystretch}{0.7}\begin{array}{l}
\framebox{1}\\
\framebox{5} \\
\end{array}$}

\put(370,32){$\framebox{4}$}

\put(326,113) {$\renewcommand{\arraystretch}{0.7}\begin{array}{l}
\framebox{1}\framebox{5}\\
\framebox{3} \\
\end{array}$}

\put(448,120) {$\renewcommand{\arraystretch}{0.7}\begin{array}{l}
\framebox{1}\framebox{5}\\
\framebox{4} \\
\end{array}$}

\put(420,110) {$\renewcommand{\arraystretch}{0.7}\begin{array}{l}
\framebox{2}\framebox{3}\\
\framebox{5} \\
\end{array}$}

\put(503,235) {$\renewcommand{\arraystretch}{0.7}\begin{array}{l}
\framebox{2}\framebox{6}\\
\framebox{5} \\
\end{array}$}

\put(470,215) {$\renewcommand{\arraystretch}{0.7}\begin{array}{l}
\framebox{2}\framebox{5}\\
\framebox{6} \\
\end{array}$}

\put(503,210) {$\renewcommand{\arraystretch}{0.7}\begin{array}{l}
\framebox{1}\framebox{7}\\
\framebox{5} \\
\end{array}$}

\put(488,245) {$\renewcommand{\arraystretch}{0.7}\begin{array}{l}
\framebox{5}\\
\framebox{7} \\
\end{array}$}

\put(313,340) {$\renewcommand{\arraystretch}{0.7}\begin{array}{l}
\framebox{2}\framebox{6}\\
\framebox{6} \\
\end{array}$}

\put(317,304) {$\renewcommand{\arraystretch}{0.7}\begin{array}{l}
\framebox{1}\framebox{7}\\
\framebox{6} \\
\end{array}$}

\put(338,333) {$\renewcommand{\arraystretch}{0.7}\begin{array}{l}
\framebox{3}\framebox{5}\\
\framebox{6} \\
\end{array}$}

\put(298,328) {$\renewcommand{\arraystretch}{0.7}\begin{array}{l}
\framebox{6}\\
\framebox{7} \\
\end{array}$}

\put(162,216) {$\renewcommand{\arraystretch}{0.7}\begin{array}{l}
\framebox{1}\framebox{6}\\
\framebox{3} \\
\end{array}$}

\put(142,230) {$\renewcommand{\arraystretch}{0.7}\begin{array}{l}
\framebox{3}\\
\framebox{6} \\
\end{array}$}

\bezier{20}(330,-70)(387,-70)(445,-70)

\bezier{20}(502,35)(473,-23)(445,-70)

\bezier{20}(502,35)(473,82)(445,130)

\bezier{20}(330,130)(387,130)(445,130)

\bezier{20}(330,130)(300,82)(273,35)

\bezier{20}(330,-70)(297,-23)(273,35)

\put(502,35){\circle{5}}

\put(445,130){\circle{5}}

\put(330,130){\circle{5}}

\put(273,35){\circle{5}}

\put(439,-82){\framebox{2}}

\put(484,34){\framebox{5}}

\put(447,140){\framebox{7}}

\put(340,138){\framebox{6}}

\put(280,32){\framebox{3}}

\bezier{120}(330,-70)(416,-22)(502,35)

\bezier{120}(502,35)(502,130)(502,230)

\bezier{120}(502,230)(416,285)(330,325)

\bezier{120}(330,325)(245,285)(163,230)

\bezier{120}(163,230)(163,130)(163,35)

\bezier{120}(163,35)(246,-13)(330,-70)

\put(502,230){\circle{5}}

\put(330,325){\circle{5}}

\put(163,230){\circle{5}}

\put(502,33) {$\renewcommand{\arraystretch}{0.7}\begin{array}{l}
\framebox{2}\\
\framebox{5} \\
\end{array}$}

\put(320,147) {$\renewcommand{\arraystretch}{0.7}\begin{array}{l}
\framebox{2}\\
\framebox{6} \\
\end{array}$}

\put(305,125) {$\renewcommand{\arraystretch}{0.7}\begin{array}{l}
\framebox{1}\\
\framebox{7} \\
\end{array}$}

\put(502,33) {$\renewcommand{\arraystretch}{0.7}\begin{array}{l}
\framebox{2}\\
\framebox{5} \\
\end{array}$}

\put(502,33) {$\renewcommand{\arraystretch}{0.7}\begin{array}{l}
\framebox{2}\\
\framebox{5} \\
\end{array}$}

\put(427,145) {$\renewcommand{\arraystretch}{0.7}\begin{array}{l}
\framebox{2}\\
\framebox{7} \\
\end{array}$}

\put(320,147) {$\renewcommand{\arraystretch}{0.7}\begin{array}{l}
\framebox{2}\\
\framebox{6} \\
\end{array}$}

\put(305,125) {$\renewcommand{\arraystretch}{0.7}\begin{array}{l}
\framebox{1}\\
\framebox{7} \\
\end{array}$}

\put(252,33) {$\renewcommand{\arraystretch}{0.7}\begin{array}{l}
\framebox{1}\\
\framebox{4} \\
\end{array}$}

\put(502,33) {$\renewcommand{\arraystretch}{0.7}\begin{array}{l}
\framebox{2}\\
\framebox{5} \\
\end{array}$}

\put(320,147) {$\renewcommand{\arraystretch}{0.7}\begin{array}{l}
\framebox{2}\\
\framebox{6} \\
\end{array}$}

\put(305,125) {$\renewcommand{\arraystretch}{0.7}\begin{array}{l}
\framebox{1}\\
\framebox{7} \\
\end{array}$}

\put(148,17) {$\renewcommand{\arraystretch}{0.7}\begin{array}{l}
\framebox{1}\\
\framebox{3} \\
\end{array}$}

\end{picture}


\end{document}